\documentclass[twoside,11pt]{article}

\usepackage{graphicx, caption, subcaption}
\usepackage[top=1in,bottom=1in,left=1in,right=1in]{geometry}
\usepackage[sort&compress, numbers]{natbib} 
\usepackage{amsfonts, amsmath, amssymb, amsthm} \usepackage{MnSymbol}
\usepackage[mathscr]{euscript}
\usepackage{mathtools}
\mathtoolsset{showonlyrefs}

\usepackage{color}
\definecolor{darkred}{RGB}{100,0,0}
\definecolor{darkgreen}{RGB}{0,100,0}
\definecolor{darkblue}{RGB}{0,0,150}

\usepackage{hyperref}
\hypersetup{colorlinks=true, linkcolor=darkred, citecolor=darkgreen, urlcolor=darkblue}
\usepackage{url}

\def\dim{v}

\newtheorem{thm}{Theorem}

\theoremstyle{remark}

\newtheorem{rem}{Remark}

\newcommand{\thmref}[1]{Theorem~\ref{thm:#1}}

\newcommand{\secref}[1]{Section~\ref{sec:#1}}
\newcommand{\figref}[1]{Figure~\ref{fig:#1}}

\def\beq{\begin{equation}} % \setcounter{equation}{1}}
\def\eeq{\end{equation}}
\def\beqn{\begin{eqnarray*}}
\def\eeqn{\end{eqnarray*}}
\def\Bitem{\begin{itemize}\setlength{\itemsep}{.2in}}
\def\bitem{\begin{itemize}\setlength{\itemsep}{.05in}}
\def\eitem{\end{itemize}}
\def\Benum{\begin{enumerate}\setlength{\itemsep}{.2in}}
\def\benum{\begin{enumerate}\setlength{\itemsep}{.05in}}
\def\eenum{\end{enumerate}}
\def\bmult{\begin{multline*}}
\def\emult{\end{multline*}}
\def\bcenter{\begin{center}}
\def\ecenter{\end{center}}
\def\bframe{\begin{frame}}
\def\eframe{\end{frame}}

\def\cA{\mathcal{A}}

\def\cC{\mathcal{C}}

\def\cE{\mathcal{E}}

\def\bbI{\mathbb{I}}

\def\bbR{\mathbb{R}}

\newcommand{\E}{\operatorname{\mathbb{E}}}
\renewcommand{\P}{\operatorname{\mathbb{P}}}

\renewcommand{\>}{\rangle}

\def\eps{\varepsilon}

\def\comp{\mathsf{c}}

\newcommand{\IND}[1]{\bbI\{ #1 \}}

\definecolor{purple}{rgb}{0.4,.1,.9}
\definecolor{darkblue}{rgb}{0,0,.8}

\newcommand{\erynew}[1]{{\leavevmode\color{black}{#1}}}

\newcommand\blfootnote[1]{%
  \begingroup
  \renewcommand\thefootnote{}\footnote{#1}%
  \addtocounter{footnote}{-1}%
  \endgroup
}

\begin{document}

\title{On the Estimation of Latent Distances Using Graph Distances}
\author{
Ery Arias-Castro\footnote{~Dept.~of Mathematics, University of California, San Diego (USA)}
\and
Antoine Channarond\footnote{~LMRS -- UMR 6085 CNRS -- Universit\'e de Rouen (France)}
\and
Bruno Pelletier\footnote{~IRMAR -- UMR CNRS 6625 -- Universit\'e Rennes II (France)}
\and
Nicolas Verzelen\footnote{~INRA -- UMR 729 MISTEA -- Montpellier (France)}
}
\date{}
\maketitle

\blfootnote{This work was partially supported by the US National Science Foundation (DMS 1513465, DMS 1916071).}

\begin{abstract}
We are given the adjacency matrix of a geometric graph and the task of recovering the latent positions.  We study one of the most popular approaches which consists in using the graph distances and derive error bounds under various assumptions on the link function.  In the simplest case where the link function is proportional to an indicator function, the bound matches an information lower bound that we derive.
\end{abstract}

\section{Introduction}

Suppose that we observe a undirected graph with adjacency matrix $W = (W_{ij} : i,j \in [n])$ (where $[n] := \{1, \dots, n\}$ and $n \ge 3$) with $W_{ij} \in \{0,1\}$ and $W_{ii} = 0$.  We assume the existence of points, $x_1, \dots, x_n \in \bbR^\dim$, such that 
\beq\label{W}
\P(W_{ij} = 1 \mid x_1, \dots, x_n) = \phi(\|x_i - x_j\|),
\eeq 
for some non-increasing link function $\phi : [0, \infty) \mapsto [0,1]$.
The $(W_{ij}, i < j)$ are assumed to be independent given the point set $(x_1, \dots, x_n)$.  
We place ourselves in a setting where the adjacency matrix $W$ is observed, but the underlying points are unknown.  
We will be mostly interested in settings where $\phi$ is unknown (and no parametric form is known).
Our most immediate interest is in the pairwise distances 
\beq\label{distances}
d_{ij} := \|x_i - x_j\|.
\eeq

In general, when the link function is unknown, all we can hope for is to rank these distances.  Indeed, the most information we can aspire to extract from $W$ is the probability matrix $P := (p_{ij})$, where
\beq\label{p}
p_{ij} := \P(W_{ij} = 1 \mid x_1, \dots, x_n),
\eeq
and even with perfect knowledge of $P$, the distances can only be known up to a monotone transformation, since $p_{ij} = \phi(d_{ij})$ and $\phi$ is in principle an arbitrary non-increasing function. 
Recovering the points based on such a ranking amounts to a problem of ordinal embedding (aka, non-metric multidimensional scaling), which has a long history \cite{young1987multidimensional, kruskal1964multidimensional, shepard1962analysisI, shepard1962analysisII}.  

Although this is true in general, we focus our attention on the `local setting' where the link function has very small support.  In that particular case, we are able to (approximately) recover the pairwise distances up to a scaling.
By fixing the scale arbitrarily (since it cannot be inferred from the available data), recovering the underlying points amounts to a problem of metric multidimensional scaling \cite{borg2013modern}. Classical Scaling \cite{torgerson1952multidimensional} is the most popular method for that problem, and comes with a perturbation bound \cite{arias2018perturbation} which can help translate an error bound for the estimation of the pairwise distances (up to scale) to an error bound for the estimation of the points (up to a similarity transformation).
We thus focus our attention on the estimation of the pairwise distances \eqref{distances}.

\subsection{Related work}
The model we consider in \eqref{W} is an example of a latent graph model and the points are often called latent positions.  In its full generality, the model includes the planted partition model popular in the area of graph partitioning.  To see this, take $r = 1$ and let $\dim$ denote the number of blocks and, with $e_k$ denoting the $k$-th canonical basis vector, set $x_i = e_k$ if $i$ belongs to block $k$.  The planted partition model is a special case of the stochastic block model of \citet{Holland1983109}.  This is also a special case of our model, as can be seen by changing $e_k$ to $z_s$ chosen so that $\phi(\|z_s - z_\ell\|) = p_{k\ell}$, where $p_{k\ell}$ denotes the connection probability between blocks $k$ and $\ell$.  
Mixed-membership stochastic block models as in \cite{anandkumar2013tensor, 
yang2013overlapping, 
airoldi2008mixed} are also special cases of latent graph models, but of a slightly different kind.
The literature on the stochastic block model is now substantial and includes results on the recovery of the underlying communities; see, e.g., 
\cite{sussman2012consistent,
lei2015consistency,
rohe2011spectral,
mossel2015consistency,
chen2016statistical,
hajek2016achieving,
abbe2017community} and references therein.

Our contribution here is of a different nature as we focus on the situation where the latent positions are well spread out in space, forming no obvious clusters.  This relates more closely to the work of \citet{hoff2002latent}.  Although their setting is more general in that additional information may be available at each position, without that additional information their approach reduces to the following logistic regression model:
\beq\label{logistic}
\log\Big( \frac{p_{ij}}{1 - p_{ij}} \Big) = -d_{ij},
\eeq
which is clearly a special case of \eqref{W} with link function the logistic function.
\citet{sarkartheoretical} consider this same model motivated by a link prediction problem where the nodes are assumed to be embedded in space with their Euclidean distances being the dissimilarity of interest.  In fact, they assume that the points are uniformly distributed in some region.  They study a method  based on the number of neighbors that a pair of nodes have in common, which is one of the main methods for link prediction \cite{liben2007link, liben2003link}.  
\citet{parthasarathy2017quest} consider a more general setting where a noisy neighborhood graph is observed: if $(x_i)$ are points in a metric space with pairwise distances $(d_{ij})$, then an adjacency matrix, $W = (W_{ij})$, is observed, where $W_{ij} = 1$ with probability $1-p$ if $d_{ij} \le r$ and with probability $q$ if $d_{ij} > r$, where $p, q \in [0,1]$ are parameters of the model.  Under fairly general conditions on the metric space and the sampling distribution, and additional conditions on $(n,r,p)$, they show that the graph distances computed based on $W$ provide, with high probability, a 2-approximation to the underlying distances in the case where $q=0$.  In the case where $q>0$, the same is true, under some conditions on $(n,r,p,q)$, if $W$ is replaced by $\tilde W = (\tilde W_{ij})$ where $\tilde W_{ij} = 1$ exactly when $N_{ij}/(N_i+N_j-N_{ij}) \ge \tau$, where $\tau$ is a carefully chosen tuning parameter and where $N_i := \#\{j : W_{ij} = 1\}$ (number of neighbors of $i$) and $N_{ij} :=  \#\{k : W_{ik} = W_{jk} = 1\}$ (number of common neighbors of $i$ and $j$).

\citet{scheinerman2010modeling} and \citet{young2007random} consider what they call a dot-product random graph model where $p_{ij} = \<x_i, x_j\>$, where it is implicitly assumed that $\<x_i, x_j\> \in [0,1]$ for all $i \ne j$.  This model is a special case of \eqref{W}, with $\phi(d) = 1 - \frac12 d^2$.
\citet{sussman2014consistent} consider recovering the latent positions in this model with full knowledge of the link function.
They devise a spectral method which consists in embedding the items $\{1, \dots, n\}$ as points in $\bbR^\dim$, with $v$ assumed known, as the row vectors of $\smash{U_{(\dim)} \Theta_{(\dim)}^{1/2}}$, where $W = U \Theta V^\top$ is the SVD of $W$, and for a matrix $A = (A_{ij})$ and an integer $s \ge 1$, $A_{(s)} = (A_{ij} : i \vee j \le s)$. 
They analyze their method in a context where the latent positions are in fact a sample from a possibly unknown distribution.
The same authors extended their work in \cite{tang2013universally} to an arbitrary link function, which may be unknown, although the focus is on a binary classification task in a setting where for each $i \in [n]$ a binary label $y_i$ is available.

\citet{alamgir2012shortest, von2013density} consider the closely related problem of recovering the latent positions in a setting where a nearest-neighbor graph is available.   They propose a method based on estimating the underlying density denoted $f$.  If $\hat f_i$ denotes the density estimate at $x_i$, a graph is defined on $[n]$ with weights $\smash{w_{ij} = (\hat f_i^{-1/\dim} + \hat f_j^{-1/\dim})/2}$, and $d_{ij}$ is estimated by the graph distance between nodes $i$ and $j$. 

\erynew{
Latent positions random graph models also play a role in the literature on rankings\footnote{ Thanks to Philippe Rigollet for pointing out this out to us.} \cite{fligner1993probability, marden1996analyzing}. 
A typical parametric model represents each player $i \in [n]$ by a number $x_i$ such that the probability that $i$ wins against $j$ in a single game is $\phi(x_i - x_j)$. Note that the link function is applied to the difference and not the absolute value of the difference. For example, the Bradley--Terry--Luce model \cite{bradley1952rank, luce2012individual} uses the logistic link function. 
Suppose that multiple games are played between multiple pairs of players. The result of that can be summarized as $(W_{ij} : i \ne j)$, where $W_{ij}$ is the number of games where $i$ prevailed over $j$.  This is the weight matrix of a directed latent positions graph where the positions are $(x_1, \dots, x_n)$.  We refer the reader to \cite{pananjady2020worst, shah2016stochastically} and references therein for theoretical results developed for such models.
}

\subsection{Our contribution}
Graph distances are well-known estimates for the Euclidean distances in the context of graph drawing \cite{kruskal1980designing, shang2003localization}, where the goal is to embed items in space based on an incomplete distance matrix.
They also appear in the literature on link prediction \cite{liben2007link, liben2003link} and are part of the method proposed in \cite{von2013density}.  
We examine the use of graph distances for the estimation of the Euclidean distances \eqref{distances}.  As we shall see, the graph distances are directly useful when the link function $\phi$ is compactly supported, which is for example the case in the context of a neighborhood graph where $\phi(d) = \IND{d \le r}$ for some connectivity radius $r > 0$.  In fact, the method is shown to achieve a minimax lower bound in this setting (under a convexity assumption).  This setting is discussed in \secref{simple}.
In \secref{general}, we extend the analysis to other (compactly supported) link functions.
We end with \secref{discussion}, where we discuss some important limitations of the method based on graph distances and consider some extensions, including localization (to avoid the convexity assumption) and the use of the number of common neighbors (to accommodate non-compact link functions). 
%We also discuss a regularization known as Maximum Variance Unfolding (MVU) \citep{weinberger2007graph}.  
%
Proofs are gathered in \secref{proofs}.

\subsection{The graph distance method}
\label{sec:preliminaries}
Given the adjacency matrix $W$, the graph distance (aka shortest-path distance) between nodes $i$ and $j$ is defined as
\beq\label{delta}
\delta_{ij} := \inf \Big\{\ell : \exists k_0, \dots, k_\ell \in [n] \text{ s.t. } k_0 = i, k_n = j, \text{ and } W(k_{s-1}, k_s) = 1, \forall s \in [\ell]\Big\},
\eeq
where $\inf \emptyset = \infty$ by convention.
Here and elsewhere, we will sometimes use the notation $W(i,j)$ for $W_{ij}$, $d(i,j)$ for $d_{ij}$, etc.

\erynew{
We propose estimating, up to a scale factor, the Euclidean distances \eqref{distances} with the graph distances \eqref{delta}.  
Indeed, since $\phi$ is assumed unknown, the scale factor cannot be recovered from the data, as is the case in ordinal embedding, for example. Therefore, estimates are necessarily up to an arbitrary scaling factor, so that the accuracy of an estimator $\hat d = (\hat d_{ij})$ for $d = (d_{ij})$ is measured according to close we can make $s \hat d$ and $d$ in some chosen way by choosing the scale $s > 0$ with (oracle) knowledge of $d$. For example, with mean squared error, this leads to quantifying the accuracy of $\hat d$ as follows
\beq
\min_{s > 0} \sum_{i < j} (s \hat d_{ij} - d_{ij})^2.
\eeq 
}

The graph distance method is the analog of the MDS-D method of \citet{kruskal1980designing} for graph drawing, which is a setting where some of the distances \eqref{distances} are known and the goal is to recover the missing distances.  Let $\cE$ denote the set of pairs $i < j$ for which $d_{ij}$ is known.  
MDS-D estimates the missing distances with the distances in the graph with node set $[n]$ and edge set $\cE$, and with edge $(i,j) \in \cE$ weighed by $d_{ij}$.  
This method was later rediscovered by \citet{shang2003localization}, who named it MDS-MAP(P), and coincides with the IsoMap procedure of \citet{Tenenbaum00ISOmap} for isometric manifold embedding.  (For more on the parallel between graph drawing and manifold embedding, see the work of \citet{chen2009local}.)

As we shall see, the graph distance method is most relevant when the positions are sufficiently dense in their convex hull, which is a limitation it shares with MDS-D. 
For $\Omega \subset \bbR^\dim$ and $x_1, \dots, x_n \in \bbR^\dim$, define 
\beq\label{eps}
\Lambda_\Omega(x_1, \dots, x_n) = \sup_{x \in \Omega} \min_{i \in [n]} \|x - x_i\|,
\eeq
which measures how dense the latent points are in $\Omega$.
We also let  $\Lambda(x_1, \dots, x_n)$ denote \eqref{eps} when $\Omega$ is the convex hull of $\{x_1, \dots, x_n\}$.

\section{Simple setting} \label{sec:simple}

In this section we focus on the simple, yet emblematic case of a neighborhood (ball) graph, that is, a setting where the link function is given by $\phi(d) = \IND{d \le r}$ for some $r > 0$.  
\erynew{
When the positions are drawn iid from some distribution, the result is what is called a random geometric graph \cite{penrose2003random}, but here we consider the positions to be deterministic. In particular, the setting is not random.
}    

We start with a performance bound for the graph distance method and then establish minimax lower bound.  
Similar results are available in \cite{bernstein2000graph, parthasarathy2017quest, chazal2013persistence}, among other places, and we only provide a proof for completeness, and also to pave the way to the more sophisticated \thmref{graph-distance-general}.

\begin{thm} \label{thm:graph-distance-simple}
Consider a set of points $x_1, \dots, x_n$ that satisfy $\Lambda(x_1, \dots, x_n) \le \eps$. Assume that $\phi(d) = \IND{d \le r}$ for some $r > 0$, and define $\hat d_{ij} = r \delta_{ij}$.    
If the connectivity radius $r$ is sufficiently larger than the density of the point set $\eps$, specifically if $\eps \le r/4$, then 
\beq\label{graph-distance-bound}
0 \le \hat d_{ij} - d_{ij} \le 4(\eps/r) d_{ij} + r, \quad \forall i, j \in [n],
\eeq
\erynew{which in particular implies that
\beq
\left[\frac1{\binom{n}2} \sum_{i < j} (\hat d_{ij} - d_{ij})^2\right]^{1/2}
\le 4 (\eps/r) \rho + r,
\eeq
where $\rho$ is the diameter of $\{x_1, \dots, x_n\}$.}
\end{thm}
%The quantity $\rho$ is the statement is the $\ell_2$ measure of diameter, and is bounded from above by the $\ell_\infty$ measure of diameter (the usual notion).
In the statement, $\hat d$ is not a true estimator in general as it relies on knowledge of $r$, which may not be available, nor be estimable, if the link function is unknown.  Nevertheless, the result says that, up to that scale parameter, the graph distances achieve a nontrivial level of accuracy. Compare with  \cite[Th 2.5]{parthasarathy2017quest}, which in the context of points on a Euclidean space as considered here says that, in a stochastic setting where the points are generated iid from some distribution supported on a convex set, $\max_{ij} (\hat d_{ij} - d_{ij})$ is bounded by $r$ in the limit where $n \to \infty$ while $r$ remains fixed.

\begin{figure}[ht]
\centering
\begin{subfigure}[t]{0.5\textwidth}
\centering
\caption{Latent positions}
\includegraphics[scale=0.4]{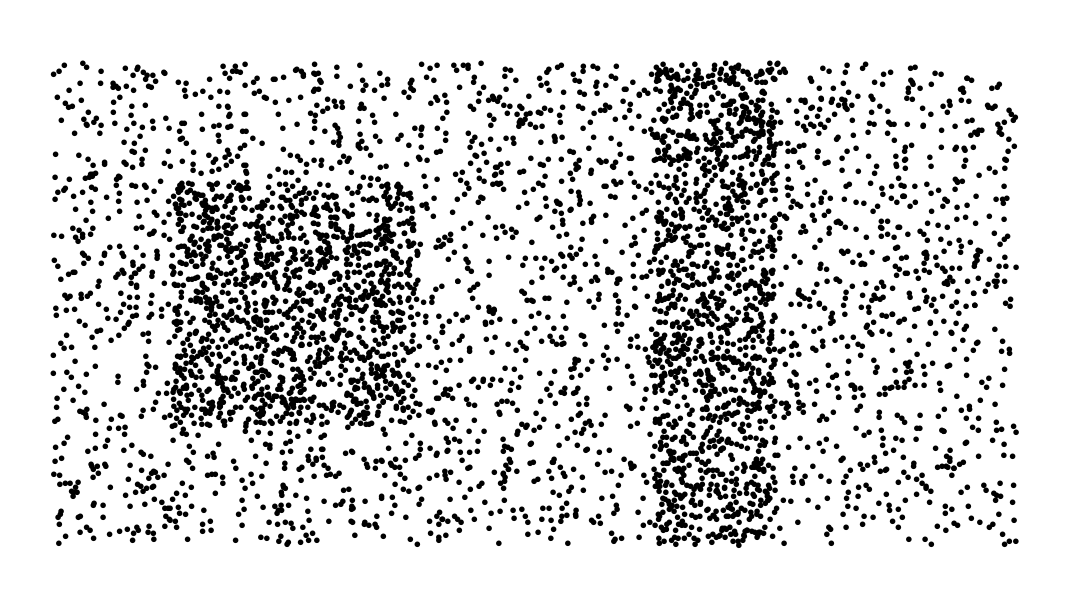}
\end{subfigure}%
\begin{subfigure}[t]{0.5\textwidth}
\centering
\caption{Recovered positions with $r = 0.05$}
\includegraphics[scale=0.4]{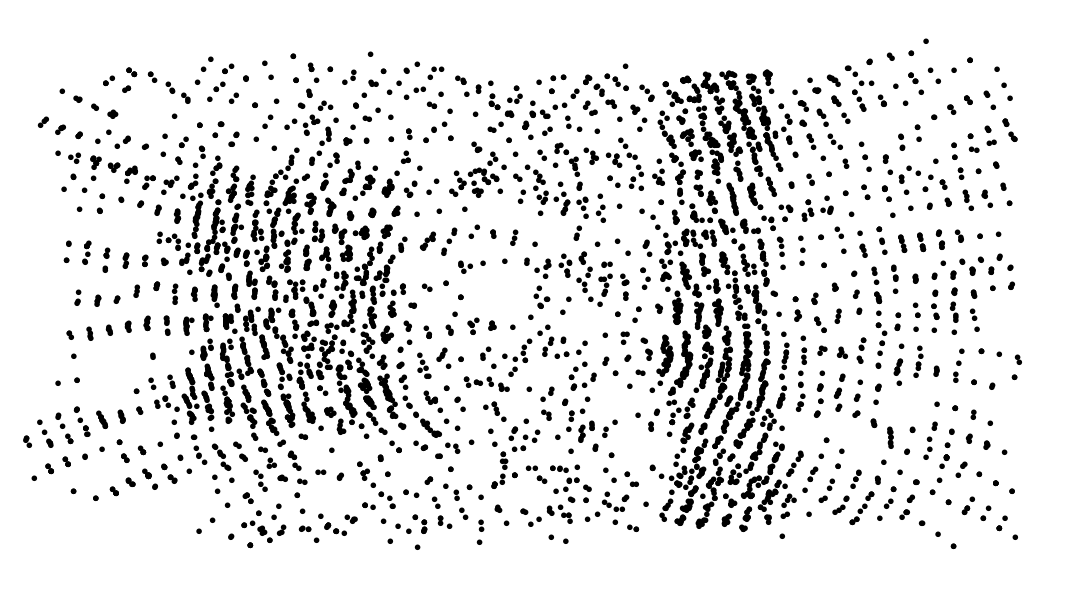}
\end{subfigure}%

\begin{subfigure}[t]{0.5\textwidth}
\centering
\includegraphics[scale=0.4]{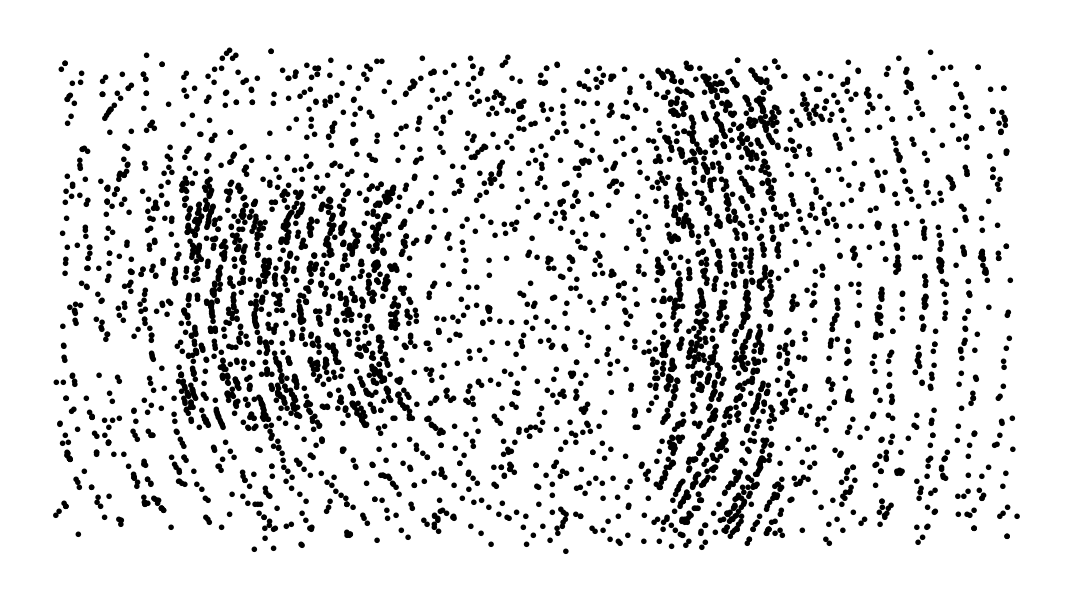}
\caption{Recovered positions with $r = 0.1$}
\end{subfigure}%
\begin{subfigure}[t]{0.5\textwidth}
\centering
\includegraphics[scale=0.4]{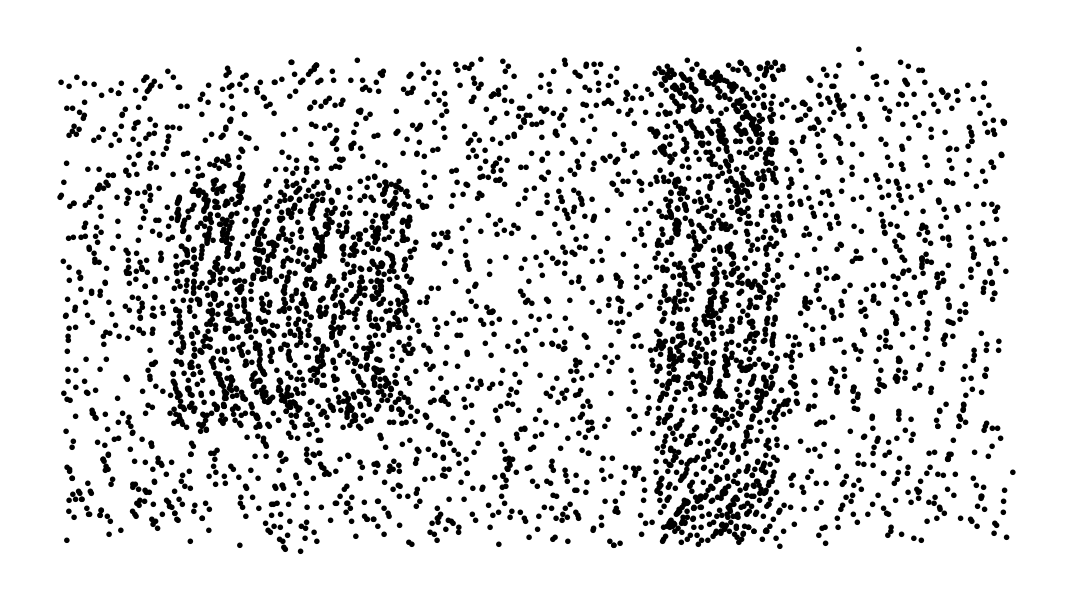}
\caption{Recovered positions with $r = 0.2$}
\end{subfigure}%
\caption{\small A numerical example illustrating the setting of \thmref{graph-distance-simple}.  Here $n_0 = 3000$ positions were sampled uniformly at random from $\Omega_0 := [0,2] \times [0,1]$, $n_1 = 1000$ from $\Omega_1 := [0.25, 0.75] \times [0.25, 0.75]$, and $n_2 = 1000$ from $\Omega_2 := [1.25, 1.5] \times [0,1]$, for a total of $n = 5000$ positions.}
\label{fig:rectangles}
\end{figure}

\begin{figure}[ht]
\centering
\begin{subfigure}[t]{0.5\textwidth}
\centering
\caption{Latent positions}
\includegraphics[scale=0.4]{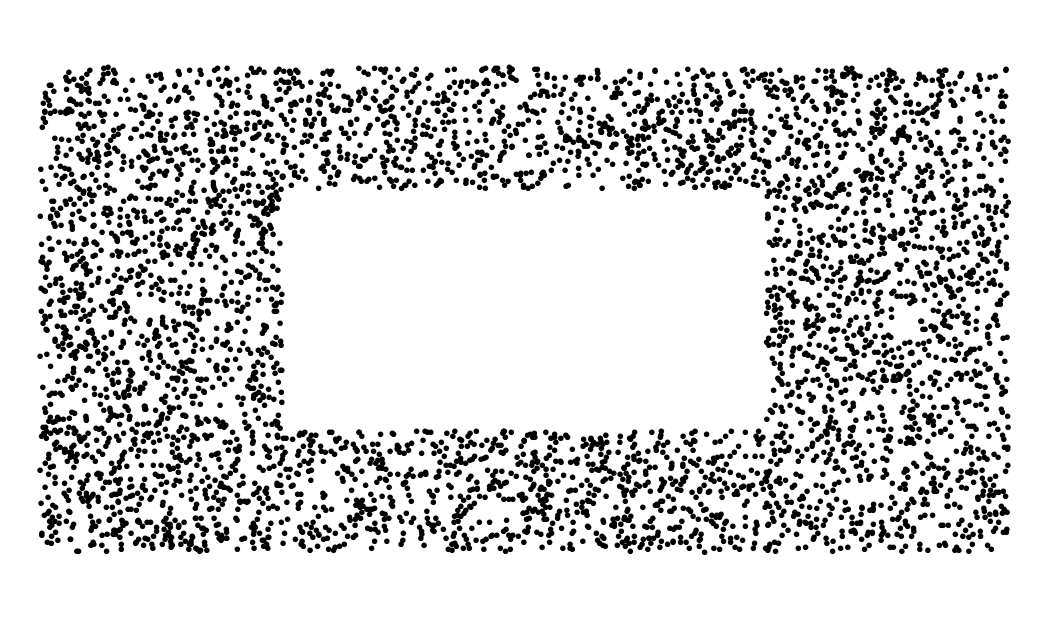}
\end{subfigure}%
\begin{subfigure}[t]{0.5\textwidth}
\centering
\caption{Recovered positions with $r = 0.2$}
\includegraphics[scale=0.4]{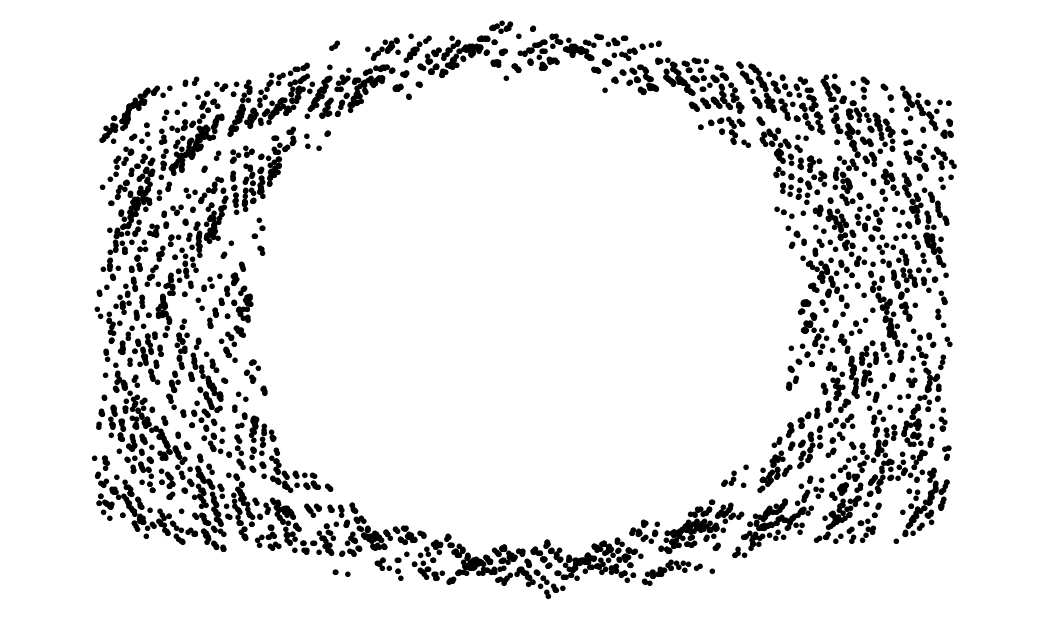}
\end{subfigure}%
\caption{\small A numerical example illustrating the setting of \thmref{graph-distance-simple} showing that the convexity constraint is indeed required for the graph distance method to be unbiased.  Here $n = 5000$ positions were sampled uniformly at random from $\Omega := [0,2] \times [0,1] \setminus [0.5,1.5] \times [0.25, 0.75]$.}
\label{fig:rectangle_hole}
\end{figure}

\begin{figure}[ht]
\centering
\begin{subfigure}[t]{0.5\textwidth}
\centering
\caption{Latent positions}
\includegraphics[scale=0.5]{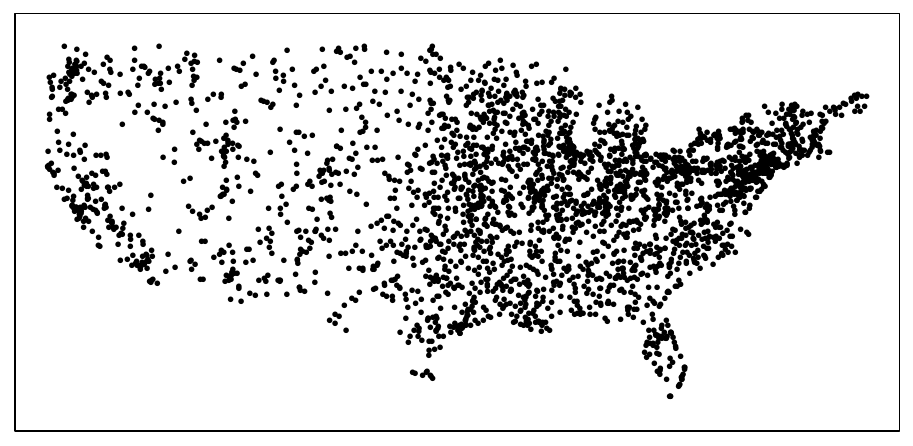}
\end{subfigure}%
\begin{subfigure}[t]{0.5\textwidth}
\centering
\caption{Recovered positions with $r = 3$}
\includegraphics[scale=0.5]{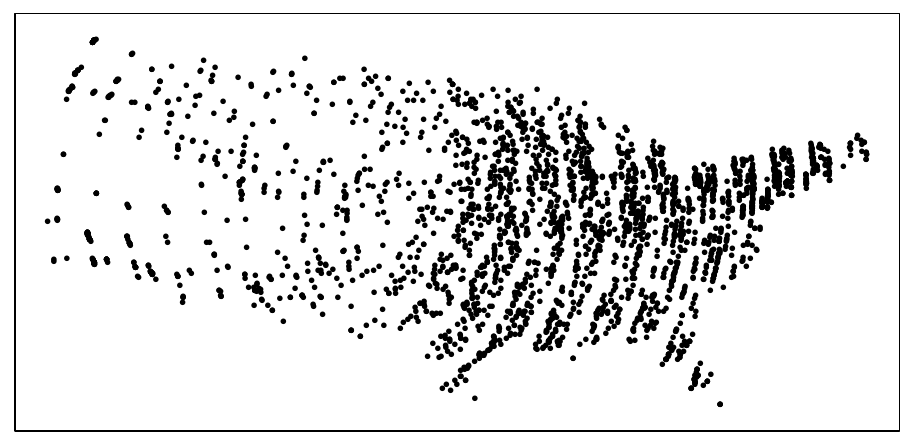}
\end{subfigure}%

\begin{subfigure}[t]{0.5\textwidth}
\centering
\includegraphics[scale=0.5]{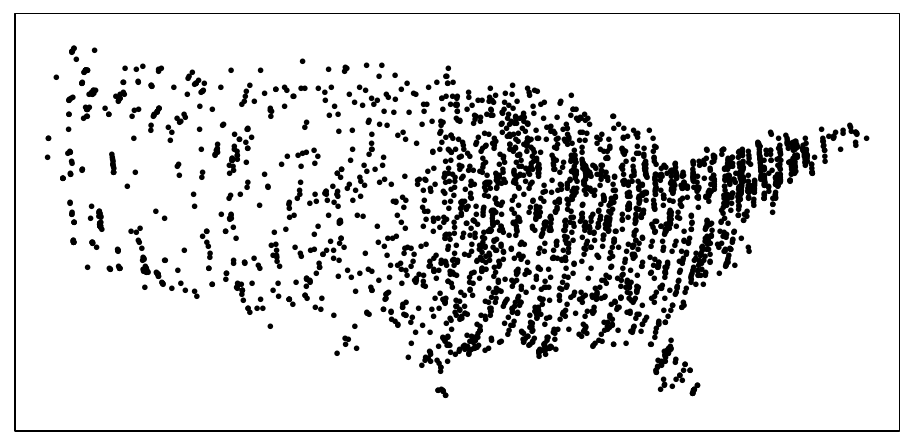}
\caption{Recovered positions with $r = 5$}
\end{subfigure}%
\begin{subfigure}[t]{0.5\textwidth}
\centering
\includegraphics[scale=0.5]{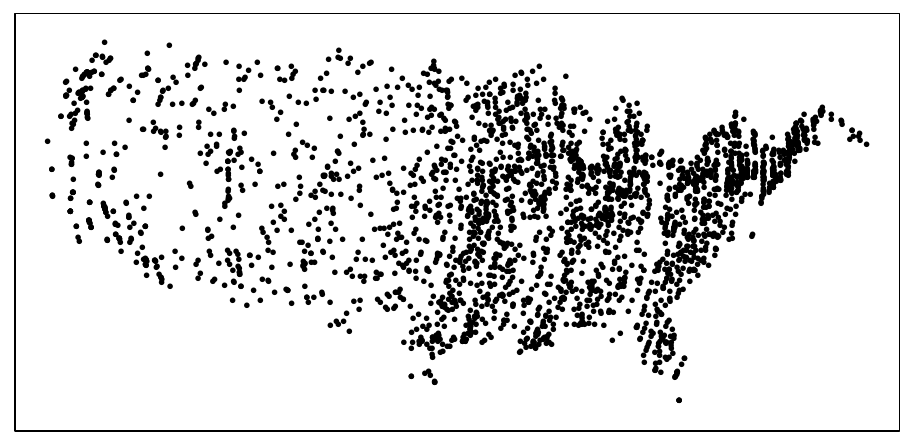}
\caption{Recovered positions with $r = 7$}
\end{subfigure}%
\caption{\small A numerical example illustrating the setting of \thmref{graph-distance-simple}.  The latent positions are located at the coordinates of $n = 3000$ US cities and the connectivity radius varies (in degrees).}
\label{fig:cities-simple}
\end{figure}

For a numerical example, see \figref{rectangles}.
In \figref{rectangle_hole} we confirm numerically that the method is biased when the underlying domain from which the positions are sampled is not convex.
That said, the method is robust to mild violations of the convex constraint, as shown in \figref{cities-simple}, where the positions correspond to $n = 3000$ US cities.\footnote{~These were sampled at random from the dataset available at \url{simplemaps.com/data/us-cities}}
(Computations were done in {\sf R}, with the graph distances computed using the {\sf igraph} package, to which Classical Scaling was applied, followed by a procrustes alignment and scaling using the package {\sf vegan}.)

\begin{rem}
If we apply Classical Scaling to $\hat d$, we obtain an embedding with arbitrary scaling and rigid positioning, which are not recoverable when $r$ is unknown.  Nevertheless, if we apply the perturbation bound recently established in \cite[Cor 2]{arias2018perturbation}, the recovery of the latent positions is of order at most $O(\eps/r + r)$. 
\end{rem}

It turns out that the graph distance method comes close to achieving the best possible performance (understood in a minimax sense) in this particularly simple setting.  Indeed, we are able to establish the following general lower bound that applies to any method.

\begin{thm} \label{thm:minimax}
\erynew{
Assume that $\phi(d) = \IND{d \le r}$ with $r > 0$ known.  
Then there is a numeric constant $c_0 > 0$  with the property that, for any $\eps > 0$ and any estimator\footnote{~An estimator here is simply a function on the set of $n$-by-$n$ symmetric binary matrices with values in $\bbR_+^{n(n-1)/2}$.} $\hat d$, there is a point set $x_1, \dots, x_n$ such that $\Lambda(x_1, \dots, x_n) \le \eps$ and, for at least half of the pairs $i \ne j$,
\beq\label{minimax_1}
|\hat d_{ij} - d_{ij}| \ge c_0 \left(\frac{\, \eps}{r \vee \eps} d_{ij} + r\right),
\eeq
and also, for another numeric constant $c_1 > 0$,
\beq
\left[\frac1{\binom{n}2} \sum_{i < j} (\hat d_{ij} - d_{ij})^2\right]^{1/2}
\ge c_1 \left(\frac{\eps}{r \vee \eps} \rho + r\right),
\eeq
where $\rho$ is the diameter of the point set.
}
\end{thm}

\erynew{Thus, in the strictest sense, the graph distance method is, for this particular link function, minimax optimal (in order of magnitude). 
It turns out that the point configurations that we consider in the proof are all embedded on the real line and thus, in principle, can be embedded in any Euclidean space. It is also the case that, for these particular configurations, it does not help if we know that $x_1 < \cdots < x_n$.}

\erynew{We speculate that a better error rate can be achieved (in probability) under a stochastic model, for example, when $x_1, \dots, x_n$ are drawn iid from the uniform distribution on some `nice' domain of a Euclidean space. On the other hand, we anticipate that our performance analysis of the graph distance method is essentially tight even then. To achieve a better performance, more sophisticated methods need to be considered. Methods using neighbors-of-neighbors information \cite{sarkartheoretical, parthasarathy2017quest} are particularly compelling, but in principle require knowing (or perhaps estimating) the underlying density if it is unknown. We probe this question a little further in \secref{common neighbors} with some simple but promising numerical experiments. (All we know about this approach is that it can lead to a 2-approximation \cite{parthasarathy2017quest}.)}

%%%%
\section{General setting} \label{sec:general}
\begin{figure}[ht]
\centering
\begin{subfigure}[t]{0.5\textwidth}
\centering
\caption{Latent positions}
\includegraphics[scale=0.5]{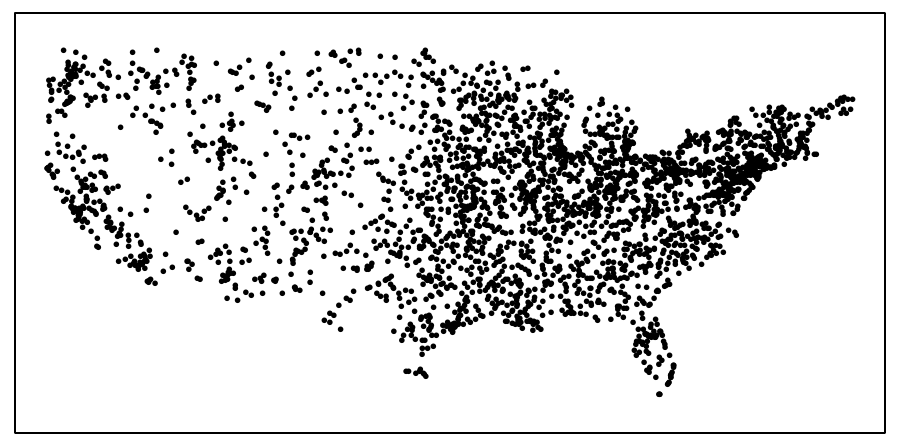}
\end{subfigure}%
\begin{subfigure}[t]{0.5\textwidth}
\centering
\caption{Recovered positions with $p = 1$}
\includegraphics[scale=0.5]{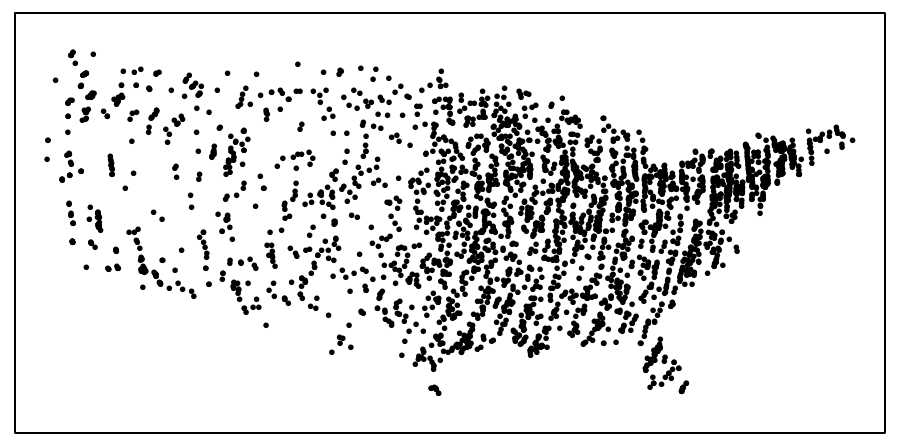}
\end{subfigure}%

\begin{subfigure}[t]{0.5\textwidth}
\centering
\includegraphics[scale=0.5]{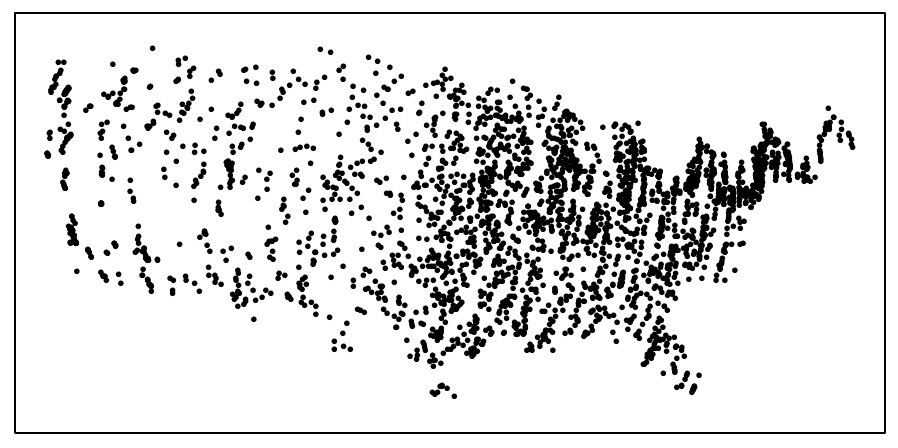}
\caption{Recovered positions with $p = 0.5$}
\end{subfigure}%
\begin{subfigure}[t]{0.5\textwidth}
\centering
\includegraphics[scale=0.5]{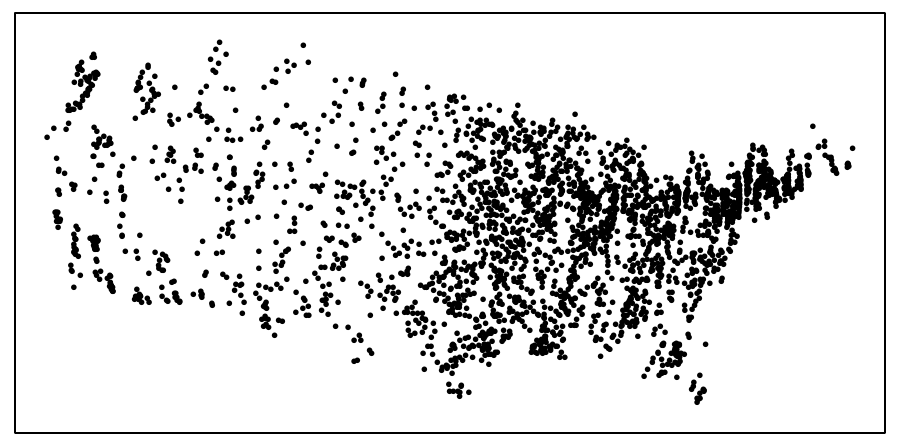}
\caption{Recovered positions with $p = 0.2$}
\end{subfigure}%
\caption{\small Same setting as in \figref{cities-simple}.  Here we set $r = 5$ and vary $p$.  In fact, to ease the comparison, we coupled the different adjacency matrices in the sense that the $(p = 0.2)$-matrix was built by erasing edges from the $(p = 0.5)$-matrix independently with probability $0.2/0.5 = 0.4$.}
\label{fig:cities-general}
\end{figure}

\erynew{
Beyond the setting of a neighborhood graph considered in \secref{simple}, the graph distance method, in fact, performs similarly when the link function is discontinuous at the edge of its support, meaning when it drops abruptly to 0.  
A case in point is when $\phi(d) = p \IND{d \le r}$ for $p > 0$ constant, which corresponds to a random geometric graph with its edges independently deleted with probability $1-p$.
See \figref{cities-general} for a numerical example illustrating this particular case.

More generally, we establish the performance of the graph distance method when the link function is compactly supported. The bound we obtain is in terms of how fast the function approaches 0 at the edge of its support.
Note that, unlike the setting of a neighborhood graph, the model is truly random when the link function is not an indicator function, so that the statement below is in probability.
}

\begin{thm} \label{thm:graph-distance-general}
Assume that $\phi$ has support $[0,r]$, for some $r > 0$, and define $\hat d_{ij} = r \delta_{ij}$.  Assume that, for some $C_0 >0$ and $\alpha \ge 0$, $\phi(d) \ge C_0 (1-d/r)^\alpha$ for all $d \in [0,r]$.  
Then there are $C_1, C_2 > 0$ depending only on $(\alpha, C_0)$ such that, whenever $r/\eps \ge C_1 (\log n)^{1+\alpha}$, for any points $x_1, \dots, x_n$ that satisfy $\Lambda(x_1, \dots, x_n) \le \eps$, with probability at least $1 - 1/n$, 
\beq\label{graph-distance-general1}
0 \le \hat d_{ij} - d_{ij} \le C_2 \big[ (\eps/r)^{\frac1{1+\alpha}} d_{ij} + r \big], \quad \forall i \ne j,
\eeq
which in particular implies that
\beq
\left[\frac1{\binom{n}2} \sum_{i < j} (\hat d_{ij} - d_{ij})^2\right]^{1/2}
\le C_2 \big[(\eps/r)^{\frac1{1+\alpha}} \rho + r\big],
\eeq
where $\rho$ is the diameter of $\{x_1, \dots, x_n\}$.
\end{thm}

\begin{figure}[ht]
\centering
\begin{subfigure}[t]{0.5\textwidth}
\centering
\caption{Median relative error versus $\alpha$.}\includegraphics[scale=0.5]{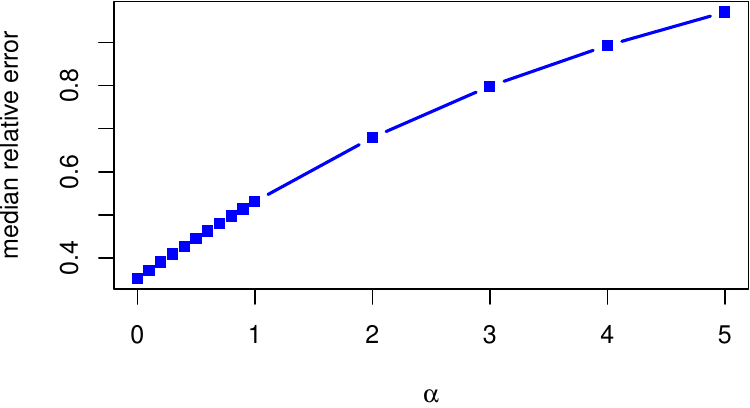}
\end{subfigure}\hfill
\begin{subfigure}[t]{0.5\textwidth}
\centering
\caption{Median relative error versus $(\eps/r)^{1/(1+\alpha)}$.}
\includegraphics[scale=0.5]{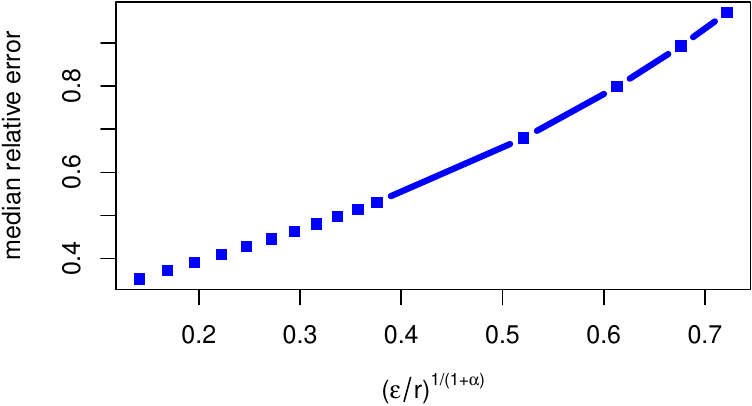}
\end{subfigure}%
\caption{\small Here $n = 5000$ points were drawn uniformly at random from the 2D unit torus. The radius was set at $r = 0.1$ and the link function varied with $\alpha$ as specified in \eqref{phi_alpha}.  The median relative error is over 100 repeats. The density $\eps$ set at $0.025$, as determined by simulation.}
\label{fig:torus}
\end{figure}

\erynew{
Although we believe our performance analysis in \eqref{graph-distance-general1} to be tight, we do not know whether it is minimax optimal in any way.  

For the graph distance, we expect it to be less accurate the slower the link function $\phi$ approaches $0$. This is borne out in some numerical experiments that we performed. In those experiments, $n = 5000$ points were drawn uniformly at random from $[0,1]^2$ considered as a torus to avoid boundary effects. (Clearly, our results apply in this setting as well.) For each $\alpha \in \{0, 0.1, \dots, 0.9, 1, 2, 3, 4, 5\}$, we computed a realization of the adjacency matrix with link function 
\begin{align}
\phi_\alpha(d) 
:= c_\alpha \big[1 \wedge (2 - 2d/r)_+^\alpha\big], \quad
c_\alpha 
:= \frac12 \frac{\alpha^2 + 3 \alpha + 2}{\alpha^2 + 5 \alpha + 8},
\end{align}
chosen so that $\P(W_{ij} = 1)$ is the same regardless of $\alpha$ as long as $r \le 0.5$. In our experiments, we chose $r = 0.1$. This was repeated 100 times. The results are presented in \figref{torus}, where 
\beq\label{phi_alpha}
\text{relative error} := \left[\sum_{i < j} (\hat d_{ij} - d_{ij})^2/\sum_{i < j} d_{ij}^2\right]^{1/2}.
\eeq
}

%\begin{figure}[ht]
%\centering
%\begin{subfigure}[t]{0.5\textwidth}
%\centering
%\caption{Latent positions}
%\includegraphics[scale=0.5]{figures/cities-x-general}
%\end{subfigure}%
%\begin{subfigure}[t]{0.5\textwidth}
%\centering
%\caption{Recovered positions with $p = 1$}
%\includegraphics[scale=0.5]{figures/cities-y-p10-r5}
%\end{subfigure}%
%
%\begin{subfigure}[t]{0.5\textwidth}
%\centering
%\includegraphics[scale=0.5]{figures/cities-y-p05-r5}
%\caption{Recovered positions with $p = 0.5$}
%\end{subfigure}%
%\begin{subfigure}[t]{0.5\textwidth}
%\centering
%\includegraphics[scale=0.5]{figures/cities-y-p02-r5}
%\caption{Recovered positions with $p = 0.2$}
%\end{subfigure}%
%\caption{\small Same setting as in \figref{cities-simple}.  Here we set $r = 5$ and vary $p$.  (In fact, to ease the comparison, we coupled the different adjacency matrices in the sense that the $(p = 0.2)$-matrix was built by erasing edges from the $(p = 0.5)$-matrix independently with probability $0.2/0.5 = 0.4$.}
%\label{fig:cities-general}
%\end{figure}

\section{Discussion} \label{sec:discussion}
The method based on graph distances suffers from a number of serious limitations:
\benum  \setlength{\itemsep}{0in}
\item The positions need to span a convex set, although the method is robust to mild violations of this constraint as exemplified in \figref{cities-simple}.
\item Even in the most favorable setting of \secref{simple}, the relative error is still of order $r$, as established in  \eqref{graph-distance-general1}. This is clearly tight for the graph distance method, and although it matches the lower bound established in \thmref{minimax}, this bias could potentially be avoided when the positions are nicely spread out, for example, as a random sample from some nice distribution is expected to be.  
\item The link function needs to be compactly supported.  
Indeed, the method can be grossly inaccurate in the presence of long edges, as in the interesting case where the link function is of the form 
\beq\label{phi2}
\phi(d) = p \IND{d \le r} + q \IND{d > r}, 
\eeq
where $0 < q < p \le 1$, as considered in \cite{parthasarathy2017quest}.
\eenum

We address each of these three issues in what follows.

\subsection{Localization}
A possible approach to addressing Issue 1 is to operate locally.  This is well-understood and is what lead  \citet{shang2004improved} to suggest MDS-MAP(P), which effectively localizes MDS-MAP \cite{shang2003localization}.  (As we discussed earlier, the latter is essentially a graph-distance method and thus bound by the convexity constraint.)
More recent methods for graph drawing based on `synchronization' also operate locally \cite{cucuringu2012sensor, cucuringu2013asap}.

\begin{figure}[ht]
\centering
\begin{subfigure}[t]{0.5\textwidth}
\centering
\caption{Latent positions}
\includegraphics[scale=0.4]{figures/rectangle_hole_n5000}
\end{subfigure}%
\begin{subfigure}[t]{0.5\textwidth}
\centering
\caption{Recovered positions with $r = 0.2$ and $\delta_{ij} \le 2$}
\includegraphics[scale=0.4]{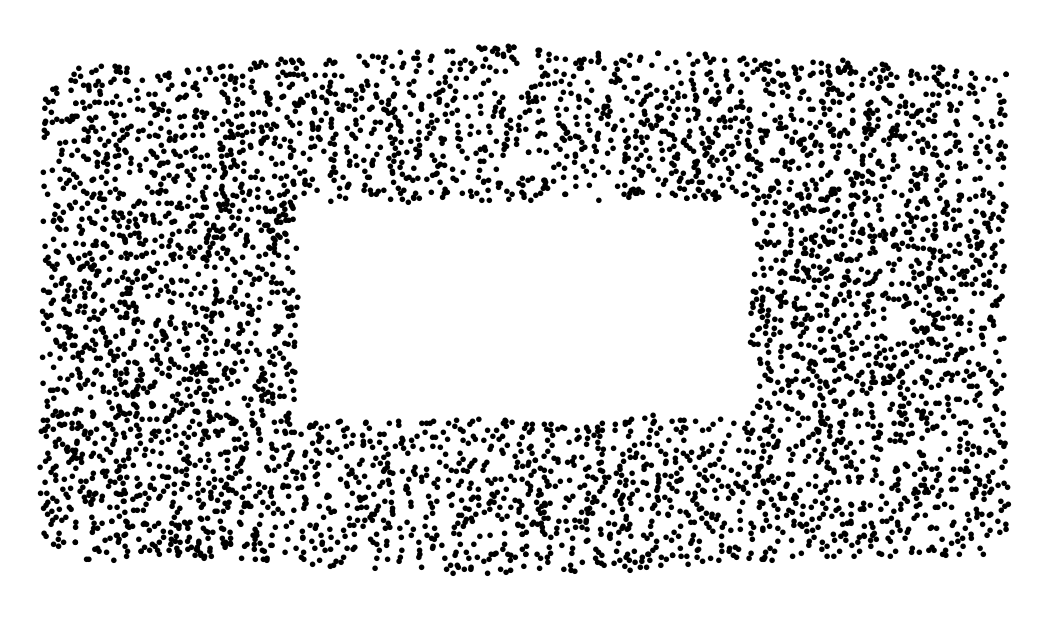}
\end{subfigure}%
\caption{\small Same setting as in \figref{rectangle_hole}.}
\label{fig:rectangle_hole_local}
\end{figure}

Experimentally, this strategy works well.  See \figref{rectangle_hole_local} for a numerical example, which takes place in the context of the rectangle with a hole of \figref{rectangle_hole}.  We adopted a simple approach: we kept the graph distances that were below a threshold, leaving the other ones unspecified, and then applied a method for multidimensional scaling with missing values, specifically SMACOF \cite{de2009multidimensional} (initialized with the output of the graph distance method).

\subsection{Regularization}

\begin{figure}[ht]
\centering
\begin{subfigure}[t]{0.5\textwidth}
\centering
\caption{Latent positions}
\includegraphics[scale=0.4]{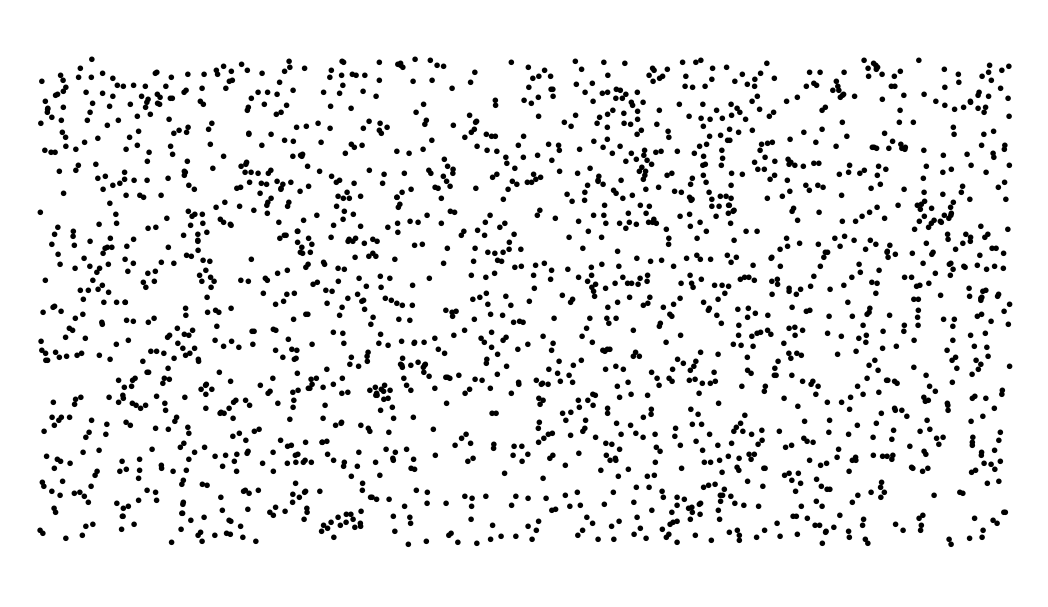}
\end{subfigure}%
\begin{subfigure}[t]{0.5\textwidth}
\centering
\caption{Recovered positions with $r = 0.5$}
\includegraphics[scale=0.4]{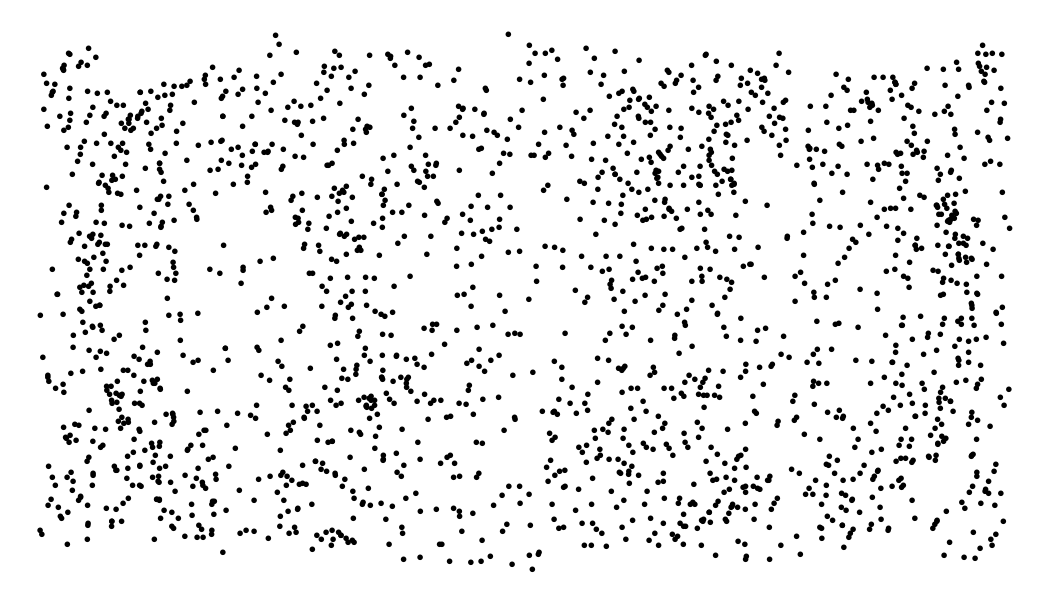}
\end{subfigure}%
\caption{\small Here $n = 2000$ positions were sampled uniformly at random from $\Omega := [0,2] \times [0,1]$.  For all $i \ne j$, $\delta_{ij} \in \{1, 2, 3, 4, 5\}$, so that the graph distances are rather discrete, yet the embedding computed by classical multidimensional scaling is surprisingly accurate.}
\label{fig:rectangle_regularization_mds}
\end{figure}

%\paragraph{Multidimensional scaling}
Regarding Issue 2, in numerical experiments we have found that the graph distances, although grossly inaccurate, are nevertheless useful for embedding the points using (classical) multidimensional scaling. Thus, if one is truly interested in estimating the Euclidean distances, one may use graph distances as rough proxies for the underlying distances, apply multidimensional scaling, and then compute the distances between the embedded points.
For a numerical illustration, see \figref{rectangle_regularization_mds}.  
This phenomenon remains surprising to us and we do not have a good understanding of the situation.

\subsection{Number of common neighbors}
\label{sec:common neighbors}

A possible approach to addressing Issue 3, as well as Issue 2, is to work with the number of common neighbors, which provides an avenue to `super-resolution' in a way, at least when the positions are sampled iid from a known distribution such as the uniform distribution on a domain (known and convex).
By this we mean that, say in the simple setting of \secref{simple}, although the adjacency matrix only tells whether two positions are within distance $r$, it is possible to gather all this information to refine this assessment.  
Similarly, in the setting where \eqref{phi2} is the link function, it is possible to tell whether two positions are nearby or not.  
This sort of concentration is well-known to the expert and seems to be at the foundation of spectral methods (see, e.g., \cite[Prop 4.2]{sussman2014consistent}). 
We refer the reader to \cite{sarkartheoretical, parthasarathy2017quest}, where such an approach is considered in greater detail.

\section{Proofs} \label{sec:proofs}

\subsection{Proof of \thmref{graph-distance-simple}}
Fix $i, j \in [n]$ distinct. 

Let $m := \lfloor d_{ij}/(r-2 \eps) \rfloor$ and note that $m (r - 2\eps) \le d_{ij} \le (m+1) (r -2\eps)$.
For $s \in \{0,\dots,m+1\}$, let $z_s = x_i + \frac{s}{m+1} (x_j - x_i)$.  
We have $z_0 = x_i$ and $z_{m+1} = x_j$, and $z_0, z_1, \dots, z_{m+1}$ are on the line joining $x_i$ and $x_j$ and satisfy $\|z_s - z_{s+1}\| \le r - 2 \eps$ for all $s$.
Let $x_{k_s}$ be such that $\|z_s - x_{k_s}\| \le \eps$, with $x_{k_0} = x_i$ and $x_{k_{m+1}} = x_j$.  
Note that $x_{k_s}$ is well-defined since $z_s$ belongs to the convex hull of $\{x_1, \dots, x_n\}$ and we have assumed that $\Lambda(x_1, \dots, x_n) \le \eps$.
By the triangle inequality, for all $s \in \{0,\dots,m\}$,
\beq \label{triangle-inequality-approximation-of-segment-by-points-of-config}
\|x_{k_s} -x_{k_{s+1}}\| \le \| x_{k_s} -z_s \| + \| z_s -z_{s+1} \| + \|z_{s+1} - x_{k_{s+1}} \| \le \eps + (r -2\eps) +\eps \le r.
\eeq
Hence, $(x_{k_0},x_{k_1},\dots,x_{k_{m+1}})$ forms a path in the graph, and as a consequence, $\delta_{ij} \le m+1$.  In turn, this implies that 
\beq\label{delta-ub}
\hat d_{ij} = r \delta_{ij} \le r m + r \le r \frac{d_{ij}}{r - 2\eps} + r \le d_{ij} + 4(\eps/r) d_{ij} + r,
\eeq
using the fact that $\eps \le r/4$.  

Resetting the notation, let $k_0 = i, k_1, \dots, k_\ell = j$ denote a shortest path joining $i$ and $j$, so that $\ell = \delta_{ij}$.  By the triangle inequality,
\beq\label{d-lb}
d_{ij} = 
\|x_{k_0} - x_{k_\ell}\| \le
\sum_{s=0}^{\ell-1} \|x_{k_s} - x_{k_{s+1}}\| \le
\ell r = r \delta_{ij} = \hat d_{ij},
\eeq
using the fact that $\|x_{k_s} - x_{k_{s+1}}\| \le r$ for all $s$.

\subsection{Proof of \thmref{minimax}}
\paragraph{First term on the RHS of \eqref{minimax_1}} 
We construct two point configurations that yield the same adjacency matrix and then measure the largest difference between the corresponding sets of pairwise distances.
Assume that $r \le 1/2$ (without loss of generality) and that $m := r (n-1)$ is an integer for convenience.
We define two configurations of points, both in $\Omega := [0,1]$ (so that $\dim=1$ here).
\bitem
\item {\em Configuration 1.}
In this configuration, 
\beq
x_i = \frac{i-1}{n-1}, \quad i \in [n].  
\eeq
Note that $\Lambda_\Omega(x_1, \dots, x_n) = 1/(2n-2)$.
\item {\em Configuration 2.}
In this configuration, 
\beq
x_i = \frac{(i-1) (1 -\eta(i-1))}{(n-1)(1 -\eta(n-1))}, \quad i \in [n].
\eeq
for some $\eta > 0$ chosen small later on.
When $\eta \le 1/(2n-3)$, which we assume, $x_i$ is increasing with $x_1 = 0$ and $x_n = 1$.
Note that 
\[
\Lambda_\Omega(x_1, \dots, x_n) = \frac{1-\eta}{(2n-2)(1 -\eta(n-1))}.
\]
\eitem
The two configurations coincide when $\eta = 0$, but we will choose $\eta > 0$ in what follows.
Under Configuration 1, the adjacency matrix $W$ is given by $W_{ij} = \IND{|i-j| \le m}$.  For the design matrix to be the same under Configuration 2, it suffices that $x_1$ have (exactly) $m$ neighbors (to the right) and that $x_n$ have (exactly) $m$ neighbors (to the left); this is because $i \mapsto x_i - x_{i-1}$ is decreasing in this configuration.  
These two conditions correspond to four equations, given by 
\beq
x_{m+1} - x_1 \le r, \quad 
x_{m+2} - x_1 > r, \quad
|x_n - x_{n-m}| \le r, \quad
x_n - x_{n-m-1} > r.
\eeq
We need only consider the first and fourth as they imply the other two.  After some simplifications, we see that the first one holds when $r \le 1$, while the fourth holds when $r \le 1 - 2/(n-1)$ and $\eta \le 1/(2n-3 + m(n-m-3))$.  Since $r \le 1/2$, $r \le 1 - 2/(n-1)$ when $n \ge 5$, and we choose $\eta = 1/(2n + m(n-m))$ for example.
Then $\Lambda_\Omega(x_1, \dots, x_n) \sim 1/2n$ in Configuration 2 (same as in Configuration 1).
We choose $n = n_\eps$ just large enough that $\Lambda_\Omega(x_1, \dots, x_n) \le \eps$ in both configurations.  In particular, $\eps \sim 1/n_\eps$ as $\eps \to 0$.  Since the result only needs to be proved for $\eps$ small, we may take $n$ as large as we need.  

Now that the two designs have the same adjacency matrix, we cannot distinguish them with the available information.  It therefore suffices to look at the difference between the pairwise distances.  
%Indeed, if $\smash{d_{ij}^{(k)}}$ denotes the distance between $x_i$ and $x_j$ in Configuration $k$, then for any estimator $\smash{\hat d_{ij}}$,
%\beq
%\max_{k \in \{1,2\}} |\hat d_{ij} - d^{(k)}_{ij}| \ge \frac12\, |d^{(1)}_{ij} - d^{(2)}_{ij}|, \quad \forall i < j.
%\eeq
Let $\smash{d_{ij}^{(k)}}$ denote the distance between $x_i$ and $x_j$ in Configuration $k$. 
For $i < j$, we have
\begin{align}
d_{ij}^{(1)} &= \frac{j-i}{n-1}, &
d_{ij}^{(2)} &= \frac{j-i}{n-1} \frac{1 - \eta (j+i-2)}{1 - \eta (n-1)}.
\end{align}
In particular,
\beq
\frac12 d_{ij}^{(1)} \le d_{ij}^{(2)} \le 2 d_{ij}^{(1)},
\eeq
by the fact that $\eta n \le 1/2$.
Also,
\beq
d_{ij}^{(1)} - d_{ij}^{(2)} = \frac{j-i}{n-1} \frac{\eta (i+j -n-1)}{1 - \eta (n-1)},
\eeq
which leads to
\begin{align}
|d_{ij}^{(1)} - d_{ij}^{(2)}| 
&= d_{ij}^{(1)} \frac{\eta |i+j -n-1|}{1 - \eta (n-1)} \\
&\ge d_{ij}^{(1)} \frac{\eta n}{1 - \eta (n-1)} \frac{|i+j -n-1|}n \\
&\ge C (d_{ij}^{(1)} \vee d_{ij}^{(2)}) \frac{\eps}{\eps \vee r} \frac{|i+j -n-1|}n,
\end{align}
for some universal constant $C > 0$, using the fact that $\eta n \le 1/2$ and $\eta n \asymp 1/(1 \vee r n) \asymp \eps/(\eps \vee r)$, the latter because $\eps \asymp 1/n$ in our construction. 
Since $|i+j -n-1| \ge n/10$ for most pairs of indices $i < j$, the following is also true
\beq
|d_{ij}^{(1)} - d_{ij}^{(2)}| \ge (C/10) \frac{\eps}{r \vee \eps} (d_{ij}^{(1)} \vee d_{ij}^{(2)}).
\eeq
To conclude, since the two configurations have the same adjacency matrix, they are indistinguishable solely based on that information, and so it must be that for any estimator $\hat d$, for most pairs $i < j$,
\begin{align}
|\hat d_{ij} - d_{ij}| 
&\ge |\hat d_{ij} - d_{ij}^{(1)}| \vee |\hat d_{ij} - d_{ij}^{(2)}| \\
&\ge \frac12 |d_{ij}^{(1)} - d_{ij}^{(2)}| \\
&\ge (C/20) \frac{\eps}{r \vee \eps} d_{ij},
\end{align}
where $d_{ij}$ denotes $d_{ij}^{(k)}$ if the true configuration is Configuration $k$.

\erynew{
\paragraph{Second term on the RHS of \eqref{minimax_1}} 
We construct again two point configurations, also on the real line, that have the same adjacency matrix. We assume that $\eps \le r$ for otherwise the first term in the RHS of \eqref{minimax_1} is of order $\asymp 1$ for most pairs of indices. In fact, we assume that $q = r/\eps$ is an integer for simplicity.

To any {\em pattern} $y_0 = 0 < y_1 < \cdots < y_m = r$ with $y_j - y_{j-1} \le 2\eps$ for all $j$, associate the point set $x_1 < \cdots < x_n$ where $x_i = y_{i\, {\rm mod}\, m} + \lfloor i/m \rfloor r$. Note that the $x$ point set is built by repeating the $y$ pattern. As can be readily seen, all these point sets have the same adjacency matrix $W_{ij} = \IND{|i-j| \le m}$ and have $\Lambda$ bounded by $\eps$.
We now consider two particular cases. 
Take $m$ even so that $m = 2q$ for some integer $q > 0$, and define the following configurations. 
\bitem
\item {\em Configuration 3.} 
Here $x_1, \dots, x_n$ is defined based on $y_j = j \eps$ for $j = 1, \dots, q-1$ and $y_j = (q-1)\eps + (j-q+1) \eta$ for $j = q, \dots, 2q$, where $\eta := \eps/(q+1)$. 
\item {\em Configuration 4.}
Here the configuration is obtained by reversing the order of the previous one, namely, $x_1, \dots, x_n$ is based on $y_j = j \eta$ for $j = 1, \dots, q$ and $y_j = (j-q) \eps$ for $j = q+1, \dots, 2q$.
\eitem

Let $\smash{d_{ij}^{(k)}}$ denote the distance between $x_i$ and $x_j$ in Configuration $k$. 
Letting $b_i := i\, {\rm mod }\, 2q$, we have
\begin{align}
d_{ij}^{(3)} &= (2q-b_i)\eta + (q-1-b_i)_+(\eps-\eta) + b_j\eta +b_j\wedge(q-1)(\eps-\eta) + r(\lfloor j/2q \rfloor - \lfloor i/2q \rfloor -1), \\
d_{ij}^{(4)} &= (2q-b_i)\eta + (2q-b_i)\wedge(q-1)(\eps-\eta) + b_j\eta +(b_j-q-1)_+(\eps-\eta) + r(\lfloor j/2q \rfloor - \lfloor i/2q \rfloor -1),
\end{align}
resulting in
\beq
|d_{ij}^{(3)} - d_{ij}^{(4)}|
= (\eps - \eta) \cdot \begin{cases}
|b_j-b_i| & \text{if } b_i,b_j\le q-1 \text{ or } b_i,b_j\ge q+1 \text{ or } b_i=b_j=q, \\
|2q-b_i-b_j| & \text{if } b_i\le q-1,b_j\ge q+1 \text{ or } b_i\ge q+1, b_j\le q-1, \\
|b_j-b_i-1| & \text{if } b_i=q, b_j\ge q+1 \text{ or } b_i\le q-1,b_j=q, \\
|b_j-b_i+1| & \text{if } b_i=q, b_j\le q-1 \text{ or } b_i\ge q+1,b_j=q,
\end{cases}
\eeq
where for a real $a$, $a_+ := \max(0, a)$.
We have $\eps - \eta = q \eps/(q+1) = r/(q+1)$, and elementary considerations confirm that for most $i < j$, the factor defined in the curly bracket above is $\ge C (q+1)$ for some universal constant $C>0$. Hence, $|d_{ij}^{(3)} - d_{ij}^{(4)}| \ge C r$ for most pairs of indices, and this then implies as before that for any estimator $\hat d$, for most pairs $i < j$,
\begin{align}
|\hat d_{ij} - d_{ij}| 
&\ge (C/2) r,
\end{align}
where $d_{ij}$ denotes $d_{ij}^{(k)}$ if the true configuration is Configuration $k$.
%\paragraph{Conclusion}
%We obtain the stated \eqref{minimax_1} by considering by taking the union of Configuration $k \in \{1,2\}$ and Configuration $\ell \in \{3,4\}$ translated by $1+r$ XXX
}

\subsection{Proof of \thmref{graph-distance-general}}
In the following, $C_0, C_1, C_2$ refer to the constants appearing in the statement of \thmref{graph-distance-general}, while $c_1, c_2, \dots$ denote positive constants that only depend on $(\alpha, C_0)$.
Since the result only needs to be proved for large $r/\eps$, we will take this quantity as large as needed.  
In what follows, we connect each node in the graph to itself.  This is only for convenience and has no impact on the validity of the resulting arguments.

As before in \eqref{d-lb}, we have $\hat d_{ij} \ge d_{ij}$ for all $i \ne j$.  Recall the definition of $p_{ij} \equiv p(i,j)$ in \eqref{p}.  Let $p_0 = \phi(r/2) > 0$ and note that $p_0 \ge C_0 (1/2)^\alpha$.

\bigskip\noindent
{\bf Special case.} {\em Suppose that $d_{kl} \le r/2$ for all $k \ne l$.}
In that case, for all $i \ne j$, $p_{ij} = \phi(d_{ij}) \ge \phi(r/2) = p_0$.  
For $(i,j,k)$ distinct, $(x_i, x_k, x_j)$ forms a path in the graph if and only if $W_{ik} W_{kj} = 1$, which happens with probability $p_{ik} p_{kj} \ge p_0^2$.  Therefore, by independence, 
\beq
\P(\delta_{ij} > 2) \le 
\P(W_{ik} W_{kj} = 0, \forall k \notin \{i,j\}) \le
(p_0^2)^{n-2}.
\eeq
Therefore, by the union bound, with probability at least $1 - n^2 p_0^{2n-4} \ge 1 - n^2 \exp(-c_1 n)$, we have $\delta_{ij} \le 2$, implying $\hat d_{ij} \le 2 r$, for all $i \ne j$.

\medskip\noindent
{\em Henceforth, we assume that}
\beq\label{d>r}
\max_{i \ne j} d_{ij} > r/2.
\eeq  
 
\bigskip\noindent
{\bf Claim 1.}
{\em By choosing $C_1$ large enough, the following event happens with probability at least $1 - 1/n^2$,} 
\beq
\cA_1 := \big\{\hat d_{ij} \le 9 d_{ij} + 2 r, \text{ for all } i, j \in [n]\big\}.
\eeq

Take $i, j$ such that $d_{ij} \le r/4$.
We first note that there is $j_*$ such that $d(i,j_*) > r/4$, for otherwise, for all $k,l\in[n]$, $d_{kl} \le d_{ki} + d_{il} \le r/4 + r/4 = r/2$, which would contradict our assumption \eqref{d>r}.

Define
\beq
z_s = x_i + s \frac{\eps (x_{j_*} - x_i)}{d(j_*, i)}, \quad s \in [m],
\eeq
where $m := \lfloor (r/4 -\eps)/\eps \rfloor$.  By construction each $z_s$ is on the line segment joining $x_i$ and $x_{j_*}$, and so belongs to the convex hull of $x_1, \dots, x_n$; hence, by the fact that $\Lambda(x_1, \dots, x_n) \le \eps$, there is $i_s \in [n]$ be such that $\|x_{i_s} - z_s\| \le \eps$.
By the triangle inequality,
\beq
d(i, i_s) 
= \|x_i - x_{i_s}\| 
\le \|x_i - z_s\| + \|z_s - x_{i_s}\| 
\le s \eps + \eps
\le m\eps + \eps
\le (r/4 -\eps) + \eps = r/4,
\eeq
and
\beq
d(i_s,j) \le d(i_s, i) + d(i,j) \le r/4 + r/4 = r/2.
\eeq
Therefore, for each $s \in [m]$, $(x_i,x_{i_s},x_j)$ forms a path with probability at least $p_0^2$.  
By independence, therefore, there is such an 
$s \in [m]$ with probability at least $1 - (1-p_0^2)^m$.

With the union bound and the fact that $m \ge r/5\eps$ when $r/\eps$ is large enough, we may conclude that, if $C_1$ is chosen large enough, the event
\beq
\cA_2 := \big\{\text{$\hat d_{ij} \le 2 r$ for all $i \ne j$ such that $d_{ij} \le r/4$}\big\},
\eeq
has probability at least $1 - 1/n^2$.
Indeed,
\beq
\P(\cA_2^\comp) 
\le n^2 (1-p_0^2)^m 
\le n^2 \exp(- c_2 (r/5\eps)) \le 1/n^2,
\eeq
eventually, when $r/\eps \ge C_1 (\log n)^{1+\alpha}$ with $C_1$ large enough.

Next, we prove that $\cA_2$ implies $\cA_1$, which will suffice to establish the claim.  For this, we consider the remaining case where $i, j$ are such that $d_{ij} > r/4$. 
Define $z_0 = x_i$ and 
\beq\label{z}
z_s = x_i + s \frac{\eps (x_{j} - x_i)}{d(j, i)}, \quad \text{for } s \in [m],
\eeq
where this time $m := \lfloor d_{ij}/\eps \rfloor$.
As before, for each $s \in [m]$, there is $i_s \in [n]$ such that $\|x_{i_s} - z_s\| \le \eps$.  We let $i_0 = i$ and $i_m = j$.  The latter is possible since $\|z_m -x_j\| \le \eps$.  

We have
\beq\label{d_tri_ub}
d(i_s,i_{s'}) = \|x_{i_s} - x_{i_{s'}}\| \le \|z_s - z_{s'}\| +2\eps = \eps |s -s'| +2\eps,
\eeq
so that, under $\cA_2$,
\beq
d(i_s,i_{s'}) \le r/4 \quad \text{when} \quad |s -s'| \le h := \lfloor r/4\eps -2\rfloor,
\eeq
implying that $\hat d(i_s,i_{s'}) \le 2r$ when $|s -s'| \le h$.
Thus, by the triangle inequality, under $\cA_2$, 
\beq
\hat d(i,j) 
= \sum_{k = 0}^{\lfloor m/h\rfloor -1} \hat d(i_{kh}, i_{(k+1)h}) + \hat d(i_{\lfloor m/h\rfloor h}, j)
\le \lfloor m/h\rfloor 2r + \hat d(i_{\lfloor m/h\rfloor h}, j).
\eeq
By the triangle inequality,
\beq
d(i_{\lfloor m/h\rfloor h}, j) 
\le \|z_{\lfloor m/h\rfloor h} -x_j\| + \eps
= \big|d_{ij} - \lfloor m/h\rfloor h \eps\big| + \eps,
\eeq
and it is not hard to verify that $0 \le d_{ij} - \lfloor m/h\rfloor h \eps \le (h+1)\eps$, so that $d(i_{\lfloor m/h\rfloor h}, j) \le (h+2) \eps \le r/4$, implying under $\cA_2$ that $\hat d(i_{\lfloor m/h\rfloor h}, j) \le 2r$.  
Hence, under $\cA_2$, 
\beq
\hat d(i,j) \le \lfloor m/h\rfloor 2r + 2r \le \big( (d_{ij}/\eps)/(r/4\eps - 3) + 1 \big) 2r \le 9 d_{ij} + 2 r,
\eeq
when $\eps/r$ is small enough.

We have thus established Claim 1.

\bigskip\noindent
Claim 1, of course, falls quite short of what is stated in the theorem, but we use it in the remainder of the proof.  
That said, the claim takes care of all pairs $(i,j)$ such that $d_{ij} \le 2r$.
Thus, for the remainder of the proof, we only need focus on $i, j$ such that $d_{ij} > 2r$.
Define $m$ and $z_0, \dots, z_m$ as before, and also the corresponding $i_s \in [n]$.

As before, \eqref{d_tri_ub} implies that 
\beq
d(i_s,i_{s'}) \le r \quad \text{when} \quad |s -s'| \le h := \lfloor r/\eps -2\rfloor,
\eeq
and in particular
\beq
d(i_s, j) \le r \quad \text{when} \quad s \ge m -h.
\eeq
(Note that we changed the definition of $h$.)
Similarly, we have
\beq\label{d_tri_lb}
d(i_s, j) = \|x_{i_s} - x_j\| \ge \|z_s - z_m\| -\eps = \eps |s -s'| -\eps,
\eeq
so that
\beq
d(i_s, j) > r \quad \text{when} \quad s \le m -h +2.
\eeq  

%\begin{rem}
%Some of the random variables that follow can be infinite.
%\end{rem}

For each $0\leq s\leq m$, define the random variable $H_s$ by
\beq
H_s = \max\big\{0\leq k\leq h \wedge (m-s):  W(i_s,i_{s+k})=1 \big\}.
\eeq
This is a maximum since we have set $W(k,k)=1$ for all $\in[n]$.
Note that $H_0,\dots,H_m$ are jointly independent random variables with support included in $\{0,\dots,h\}$.
Set $S_0=0$, and for $t \ge 1$, define recursively $S_t = S_{t-1} + H_{S_{t-1}}$.
Importantly, if $H_{S_t}=0$, then $S_{t'} = S_t$ for all $t'\geq t$.
Based on $\{S_t\}$, define 
\beq
T = \inf \{t: S_t > m-h\},
\eeq
with the convention that $\inf \emptyset = \infty$.

Our first objective is to bound $T$ in probability.
Given $t\geq 1$, we have
\begin{align}
\P(T>t) 
&= \P(S_t\le m-h) \\
&= \P(S_t \le m-h, H_{S_{t-1}}=0) + \P(S_t\le m-h,H_{S_{t-1}}>0).
\end{align}
On the one hand,
\begin{align}
\P(S_t\leq m-h,H_{S_{t-1}} = 0) 
&= \P(S_{t-1} \leq m-h,H_{S_{t-1}}=0) \\
&\leq \P\big(\exists\, 0 \leq s \leq m-h\,:\,H_{s}=0\big) \\
& \leq m \max_{0 \le s \le m-h} \P(H_s=0).
\end{align}
On the other hand, we note that, when $H_{S_{t-1}} > 0$, we necessarily have $S_t\geq t$, so that
\beq 
\P(S_t\leq m-h,H_{S_{t-1}}>0) 
\leq \P(t\leq S_t\leq m-h).
\eeq
Thus, 
\beq
\label{eq:ptt1}
\P(T>t) \leq m \max_{s \in [m-h]} \P(H_s=0) +  \P(t\leq S_t\leq m-h).
\eeq
In what follows, we bound each term in the right hand side of \eqref{eq:ptt1}.

\bigskip\noindent
{\bf Claim 2.1.}
{\em For any $s \in [m-h]$,
\beq
\P(H_s\leq a) \leq \P(\bar{H}\leq a), \quad a \ge 0,
\eeq
where $\bar{H}$ is a random variable supported on $\{0,1,\dots,h\}$ with distribution function
\beq
\P(\bar{H}\leq a) = \prod_{k=0}^{h-a-1} \big(1 - \phi(\epsilon(h-k)+2\epsilon)\big),
\quad 0\leq a\leq h-1.
\eeq
}
Indeed, by independence,
\begin{align}
\P(H_s\leq a) 
&= \P\big(W(i_s,i_{s+a+1})=\dots=W(i_s,i_{s+h})=0\big) \\ 
&= \prod_{k=0}^{h-a-1}\left(1 - p(i_s,i_{s+h-k})\right),
\end{align}
and by the fact that $\phi$ is non-increasing and \eqref{d_tri_ub},
\beq
p(i_s,i_{s'}) 
\geq \phi(\epsilon|s-s'| + 2\epsilon).
\eeq

\bigskip\noindent
{\bf Claim 2.2.}
{\em
We have
\beq
\label{eq:phbara}
\P(\bar{H}\leq a) \leq \exp\left( -c_3 (\eps/r)^\alpha (r/\eps-3-a)^{1+\alpha}\right), \quad 0\leq a\leq h-1.
\eeq
}

For any $0\leq a \leq h-1$, we have
\beq
\P(\bar{H}\leq a) 
\leq \exp\left( -\sum_{k=0}^{h-a-1}\phi(\epsilon(h-k)+2\epsilon)\right).
\eeq
Since $\phi$ is non-increasing, and using the lower bound we assume in the statement of the theorem,
\begin{align}
\sum_{k=0}^{h-a-1}\phi(\epsilon(h-k)+2\epsilon) 
& \geq \sum_{k=0}^{h-a-1} \int_{h-k}^{h-k+1}\phi(\epsilon y + 2\epsilon){\rm d}y\\
& = \int_{a+1}^{h+1} \phi(\epsilon y + 2\epsilon){\rm d}y\\ 
& \geq \int_{a+1}^{h+1} C_0 \left(1 - \frac{\epsilon y + 2\epsilon}{r}\right)_+^\alpha {\rm d}y\\
& = C_0 (\eps/r)^\alpha \frac{1}{1+\alpha}\left(r/\eps-3 - a\right)^{1+\alpha}.
\end{align}
This proves \eqref{eq:phbara}.

\bigskip
With Claims 2.1 and 2.2, the first term on the right-hand side of \eqref{eq:ptt1} is bounded by
\beq
m \P(\bar H = 0)
\le m \exp\left( -c_3 (\eps/r)^\alpha (r/\eps-3)^{1+\alpha}\right)
\le m \exp\left( -c_4 r/\eps\right),
\eeq
when $r/\eps$ is large enough.

\bigskip\noindent
{\bf Claim 2.3.}
{\em For any $1\leq t\leq m-h$, and any $t\leq a\le m-h$, we have
\beq
\P\left(t\leq S_t\leq a\right) \leq \P\left(\bar{H}_0+\dots+\bar{H}_{t-1} \leq a\right),
\eeq
where $\{\bar{H}_t\}$ are iid copies of $\bar{H}$.
}

First, if $t=1$, we have
\beq
\P(1\leq S_1\leq a) = \P(1\leq H_0\leq a) \leq \P(H_0\leq a) \leq \P(\bar{H}_0\leq a).
\eeq
Next, fix $t\geq 2$ and suppose that the claim is true at $t-1$.
Since $t\leq S_t\leq a$ implies that $t-1\leq S_{t-1}\leq a$, we have
\begin{align}
\P(t\leq S_t\leq a) 
&= \P(t\leq S_t\leq a, t-1\leq S_{t-1}\leq a)\\
&= \P(t\leq S_{t-1} + H_{S_{t-1}} \leq a, t-1\leq S_{t-1}\leq a)\\
&= \sum_{k=t-1}^{a} \P(t\leq k+H_k\leq a , S_{t-1}=k) \\
&\le \sum_{k=t-1}^{a} \P(H_k\leq a-k) \P(S_{t-1}=k)\\
&\le \sum_{k=t-1}^{a} \P(\bar{H}_k\leq a-k) \P(S_{t-1}=k) \\
&\le \sum_{k=t-1}^{a} \P(\bar{H}_{t-1}\leq a-k) \P(S_{t-1}=k), \\
\end{align}
where we used the fact that $\{S_{t-1}=k\}$ is independent of $H_k$, and also the fact that the $\bar{H}_t$ are iid. We may assume, and we do so, that the $\{\bar H_t\}$ are defined on the same probability space as the $\{S_t\}$, and are independent of them. 
Then,
\begin{align}
\sum_{k=t-1}^{a} \P(\bar{H}_{t-1}\leq a-k) \P(S_{t-1}=k) 
& = \P(S_{t-1}+\bar{H}_{t-1} \leq a, t-1 \leq S_{t-1}\leq a)\\
& = \sum_{k=0}^h \P(S_{t-1}+k\leq a, t-1 \leq S_{t-1}\leq a, \bar{H}_{t-1}=k) \\
& = \sum_{k=0}^h \P(t-1\leq S_{t-1}\leq a-k)\P(\bar{H}_{t-1}=k)\\
& \leq \sum_{k=0}^h \P(\bar{H}_0+\dots+ \bar{H}_{t-2}\leq a-k) \P(\bar{H}_{t-1}=k)\\
& = \P(\bar{H}_0+\dots+\bar{H}_{t-1}\leq a),
\end{align}
where the inequality comes from the recursion hypothesis.
Thus the recursion proceeds, and the claim is proved.

\bigskip\noindent
{\bf Claim 2.4.}
{\em There exists a constant $b_0>0$ and a constant $c_5>0$ such that for any $1\leq t\leq m-h$ and any $a \geq \nu b_0$, we have
\beq
\P\left(\bar{H}_0+\dots+\bar{H}_{t-1} \leq th - ta\right)
\leq
\exp\left(-c_5 t\nu^{-1}a\right),
\eeq
where $\nu := (r/\eps)^{\alpha/(\alpha+1)}$.
}

Define $\bar{U} = h - \bar{H}$ and $\bar U_t = h - \bar H_t$.  We have
\beq
\P\left(\bar{H}_0+\dots+\bar{H}_{t-1} \leq th - ta\right)
= \P\left(\bar{U}_0+\dots+\bar{U}_{t-1} \geq ta\right)
\eeq
For any $1\leq u\leq h$, we have
\beq
\P\left(\bar{U}\geq u\right) 
= \P\left(\bar{H}\leq h-u\right) 
\leq \exp\left( -c_3 (\eps/r)^\alpha (u -1)^{1+\alpha}\right),
\eeq
using \eqref{eq:phbara} and the fact that $h \le r/\eps -2$.
Hence, $\bar U$ is stochastically bounded by $1 + \nu Y$, where $Y$ is a random variable with distribution
\beq
\P(Y \ge y) = \exp(-c_3 y^{1+\alpha}), \quad y \ge 0.
\eeq
Let $\{Y_t : t \ge 1\}$ be iid with distribution that of $Y$. By construction, $\bar{U}_0+\cdots+\bar{U}_{t-1}$ is stochastically bounded $t + \nu (Y_0 + \cdots + Y_{t-1})$, implying that
\beq
\P\left(\bar{H}_0+\dots+\bar{H}_{t-1} \leq th - ta\right)
\le \P\left(Y_0 + \cdots Y_{t-1} \ge t (a-1)/\nu\right).
\eeq 
By Chernoff's bound,
\beq
\P(Y_1 + \cdots + Y_{t-1} > t b)
\le \exp(- t \zeta(b)), \quad \forall b \ge 0,
\eeq
where $\zeta$ is the rate function of $Y$. 
In particular, there is $c > 0$ such that $\zeta(b) \ge c b$ for all $b \ge b_0 := \E(Y) + 1$.  (Note that $b_0$ is just another constant.)
Thus, for $a$ such that $(a-1)/\nu \ge b_0$,
\beq
\P\left(Y_0 + \cdots Y_{t-1} \ge t (a-1)/\nu\right)
\le \exp(-t c (a-1)/\nu),
\eeq
which together with $a-1 \ge a/2$ (when $r/\eps$ is large) proves the claim.

\bigskip
With Claims 2.3 and 2.4, by choosing $a = h - (m-h)/t$, we obtain that the second term on the right-hand side of \eqref{eq:ptt1} is bounded by
\beq
\exp\left( -c_5 t \nu^{-1} a \right),
\eeq
whenever $a \ge \nu b_0$, which happens when $t \ge (m-h)/(h-\nu b_0)$. This is the case when $t \ge t^* := 2 d_{ij}/r$, which may be seen using the fact that $m-h \le m \le d_{ij}/\eps$, and that $h-\nu b_0 \ge r/\eps -3 (r/\eps)^{\alpha/(\alpha+1)} b_0 \ge r/2\eps$ when $r/\eps$ is large enough.  In fact, when $t \ge t^*$, the corresponding $a$ satisfies $a \ge r/\eps - 3 - (d_{ij}/\eps)/(2 d_{ij}/r) \ge r/2\eps$ when $r/\eps$ is large enough, so that the right-hand side of \eqref{eq:ptt1} is bounded by
\beq
\exp\left( -c_5 t \nu^{-1} (r/2\eps) \right)
\le \exp\big(- t c_6 (r/\eps)^{\frac1{1+\alpha}}\big),
\eeq
when $t \ge t^*$.

\bigskip
All combined, we have
\beq
\P(T > t) 
\le n \exp(-c_4 r/\eps) + \exp\big(- t c_6 (r/\eps)^{\frac1{1+\alpha}}\big), \quad \forall t \ge t^*,
\eeq
using the fact that $m \le n$.
In particular, this implies that 
\beq
\P(T = \infty)
\le n \exp(-c_4 r/\eps).
\eeq

Let $t^\circ = t^* + (5/c_5 b_0)\log n$.
Let $\cC$ be the probability event defined by
\beq
\cC = \left\{S_T > T h - t^\circ \nu b_0 \text{ and } T\leq t^\circ\right\}.
\eeq
We have
\begin{align}
\P\left(\cC^c\right) 
&= \P\left(S_T \le T h - t^\circ \nu b_0, T\leq t^\circ\right) + \P\left(T>t^\circ\right) \\
&\le \P\left(t^\circ \leq S_{t^\circ} \le t^\circ h - t^\circ \nu b_0\right) + \P\left(T>t^\circ\right) \\
&\le \exp(-c_5 t^\circ b_0) + n \exp(-c_4 r/\eps) + \exp\big(- t^\circ c_6 (r/\eps)^{\frac1{1+\alpha}}\big), 
\end{align}
using in the second line the fact that $S_t \ge t$ when $T \le t$, and that $S_t -th - (S_{t-1} - (t-1)h) = H_{S_{t-1}} - h \le 0$ for all $t$.
The application of Claim 2.4 here is valid when $t^\circ \le m-h$.  This is the case eventually as $t^\circ = 2 d_{ij}/r + (5/c_5 b_0) \log n$ and $m-h \ge d_{ij}/\eps - r/\eps + 1 \ge d_{ij}/3\eps$, with
\beq
\frac{2 d_{ij}/r + (5/c_5 b_0) \log n}{d_{ij}/3\eps} 
\le \eps/r + c_7 b_0 (\eps/r) \log n
\le 1,
\eeq
using the fact that $d_{ij} \ge r$ and our lower bound on $r/\eps$.
When $C_1$ is large enough, $\cC$ holds with probability at least $1-1/n^4$, eventually. 

A joint control on $V_T$ and $T$ is useful because of the following.  
Under $\{T < \infty\}$, $(i, i_{S_1}, \dots, i_{S_T})$ forms a path in the graph, so that $\hat d(i, i_{S_T}) \le T r$.  
We also have
\beq\label{d_ij lb}
d_{ij} 
\ge \sum_{t=1}^T \|z_{S_t} - z_{S_{t-1}}\|
= \sum_{t=1}^T \eps (S_t - S_{t-1})
= \eps (S_T -S_0) = \eps S_T.
\eeq
Thus, under $\cC$, 
\beq
d_{ij} 
\ge \eps S_T 
\ge \eps h T - \eps b_0 \nu t^\circ
\ge r T - (3+b_0 \nu) \eps t^\circ
\ge r T - 2 b_0 \nu \eps t^\circ,
\eeq
eventually, using the fact that $h \ge r/\eps -3$. 
Furthermore, 
\begin{align*}
\nu \eps t^\circ 
&= (r/\eps)^{\alpha/(1+\alpha)} \eps (2 d_{ij}/r + (5/c_5 b_0) \log n) \\
&= 2 (\eps/r)^{1/(1+\alpha)} d_{ij} + \big[(5/c_5 b_0) (\eps/r)^{1/(1+\alpha)} (\log n)\big] r \\
&\le 2 (\eps/r)^{1/(1+\alpha)} d_{ij} + (5/c_5 b_0) C_1^{-1/(1+\alpha)} r \\
&\le 2 (\eps/r)^{1/(1+\alpha)} d_{ij} + r,
\end{align*}
when $C_1$ is large enough.
Thus, under $\cC$, 
\beq
\hat d(i, i_{S_T}) 
\le T r 
\le d_{ij} + 2 b_0 \nu \eps t^\circ 
\le d_{ij} + c_8 ((\eps/r)^{1/(1+\alpha)} d_{ij} + r).
\eeq
And since, as in \eqref{d_tri_ub}, 
\beq
d(i_{S_T}, j)
\le \eps (m - S_T) + \eps
\le \eps (h-3) + \eps
\le \eps (r/\eps -5) + \eps
\le r.
\eeq

\bigskip\noindent
Using the union bound over all pairs suitable $(i,j)$, we have established the following.

\medskip\noindent
{\bf Claim 2.}
{\em By choosing $C_1$ large enough, with probability at least $1 - 1/n^2$, for all $i, j \in [n]$ such that $d_{ij} > 2 r$, there is $k \in [n]$ such that $\hat d_{ik} \le d_{ij} + c_8 ((\eps/r)^{1/(1+\alpha)} d_{ij} + r)$ and $d_{jk} \le r$.  (We denote this event by $\cA_3$.)}

\bigskip
Assume that $\cA_1 \cap \cA_3$ holds, which happens with probability at least $1 -2/n^2$ based on the two claims that we have established above.  In that case, let $i, j, k \in [n]$ be as in the last claim.  Because $\cA_1$ holds, $d_{jk} \le r$ implies that $\hat d_{jk} \le 11 r$.  Thus, with the union bound,
\beq
\hat d_{ij} 
\le \hat d_{ik} + \hat d_{jk}
\le d_{ij} + c_8 ((\eps/r)^{1/(1+\alpha)} d_{ij} + r) + 11r,
\eeq
so that
\beq
\hat d_{ij} - d_{ij}
\le c_9 ((\eps/r)^{1/(1+\alpha)} d_{ij} + r).
\eeq
This proves that \eqref{graph-distance-general1} holds for all $i,j \in [n]$ such that $d_{ij} > 2r$.

%\subsection{Proof of \prpref{mvu1}}
%%\begin{proof}
%We only need to consider a pair of nodes $i \ne j$ such that $\delta_{ij} < \infty$, for otherwise the bound holds by convention.  
%Let $k_0 = i, k_1, \dots, k_\ell = j$ denote a shortest path in the graph connecting $i$ and $j$, so that $\ell = \delta_{ij}$.
%We then derive
%\beq
%\gamma^*_{ij} 
%= \|y^*_i - y^*_j\|
%\le \sum_{s = 0}^{\ell-1} \|y^*_{k_{s-1}} - y^*_{k_s}\|
%\le \ell = \delta_{ij},
%\eeq
%using the triangle inequality and then the fact that $\|y^*_{k_{s-1}} - y^*_{k_s}\| \le 1$ for all $s$, since $y^*_1, \dots, y^*_n$ satisfies \eqref{mvu} and $W(k_{s-1}, k_s) = 1$ for all $s$.
%%\end{proof}
%

%\subsection*{Acknowledgments}
%The computations in this paper were done in {\sf R}, with the graph distances computed using the {\sf igraph} package, to which classical scaling was applied, followed by a procrustes alignment and scaling using the package {\sf vegan}.

\small
\bibliographystyle{abbrvnat}
\bibliography{latent-graphs}

\end{document}